\documentclass[12pt]{amsart}
\usepackage{graphics}
\input BoxedEPS.tex
 \SetTexturesEPSFSpecial
  \HideDisplacementBoxes   

\newtheorem{theorem}{Theorem}[section]
\newtheorem{lemma}[theorem]{Lemma}
\newtheorem{proposition}[theorem]{Proposition}

\newtheorem{corollary}[theorem]{Corollary}

\theoremstyle{definition}
\newtheorem{definition}[theorem]{Definition}
\newtheorem{notation}[theorem]{Notation}

\theoremstyle{remark}
\newtheorem{remark}[theorem]{Remark}

\newcommand{\NN}{ {\mathbb N} }
\newcommand{\RR}{ {\mathbb R} }
\newcommand{\CC}{{\mathbb C}}
\newcommand{\cC}{{\mathcal D}}
\newcommand{\cP}{ {\mathcal P} }
\newcommand{\cF}{{\mathcal F}}
\newcommand{\HH}{{\mathcal H}}
\newcommand{\KK}{{\mathcal{K}}}

\newcommand{\cI}{{T}}
\newcommand{\ff}{\varphi}
\newcommand{\la}{\langle}
\newcommand{\ra}{\rangle}
\newcommand{\kk}{{\text k}}
\newcommand{\tr}{\text{tr}}
\newcommand{\Tr}{\text{Tr}}
\newcommand{\EE}{\text{E}}

\newcommand{\cA}{\mathcal{A}}

\newcommand{\odo}{\otimes\cdots\otimes}
\newcommand{\cyc}{\mathbf{c}\,}
\newcommand{\SNC}{S_{NC}}
\newcommand{\NC}{NC}
\newcommand{\cyclic}{{\text{\rm cyc}}}
\newcommand{\lin}{{\text{\rm lin}}}
\def\ab{\allowbreak}
\newcommand{\oh}{\text{\rm o}}
\newcommand{\csi}{\check\psi}
\newcommand{\nappa}{\rho}

\newcounter{figurenumber}
\newcommand{\figno}{\arabic{figurenumber}%
\addtocounter{figurenumber}{1}}

\def\qed{{\unskip\nobreak\hfil\penalty50
 \hskip2em\hbox{}\nobreak\hfil 
\vbox to 7.7pt{\parindent0pt
\hsize7.7pt\hrule\vrule height7.3pt%
\hfill\vrule height7.3pt\hrule}
\parfillskip=0pt \finalhyphendemerits=0 \par}}

\title[Second order freeness]
{Second Order Freeness \\
and Fluctuations of Random Matrices: \\
I. Gaussian and Wishart matrices and cyclic Fock spaces}
\author[J. A. Mingo]{James A. Mingo $^{(*)}$}
\address{Queen's University, Department of Mathematics and
  Statistics, Jeffery Hall, Kingston, ON K7L 3N6, Canada}

\email{mingo@mast.queensu.ca}

\thanks{$^*$ Research supported by Discovery Grants and a Leadership
Support Initiative Award from the Natural Sciences and Engineering
Research Council of Canada}

\author[R. Speicher]{Roland Speicher $^{(*)(\dagger)}$}
\email{speicher@mast.queensu.ca}

\thanks{$^\dagger$ Research supported by a Premier's
  Research Excellence Award from the Province of Ontario}

\date{}

\begin{document}

\begin{abstract}

We extend the relation between random matrices and free
probability theory from the level of expectations to the
level of fluctuations.  We introduce the concept of ``second
order freeness" and interpret the global fluctuations of
Gaussian and Wishart random matrices by a general limit
theorem for second order freeness.  By introducing cyclic
Fock space, we also give an operator algebraic model for the
fluctuations of our random matrices in terms of the usual
creation, annihilation, and preservation operators. We show
that orthogonal families of Gaussian and Wishart random
matrices are asymptotically free of second order.

\end{abstract}

\maketitle

\section{Introduction}

Free probability has at least three basic facets: operator
algebras, random matrices, and the combinatorics of
non-crossing diagrams. This can be seen very clearly in
Voiculescu's generalization of Wigner's semi-circle law to
the case of several independent matrices \cite{Voi1}.  The distribution
arising in this limit of random matrices can be modelled by
a sum of creation and annihilation operators on full Fock
spaces and described nicely in terms of non-crossing
partitions.

On the random matrix side, Voiculescu's theorem describes
the leading contribution to the large $N$-limit of
expectations of traces of Gaussian random matrices. However,
in the random matrix literature there are many
investigations on more refined questions in this context. On
one side, subleading contributions to the large $N$-limit
are of interest and have to be understood up to some point
for dealing with questions concerning the largest eigenvalue
of such random matrices.  On the other side, there has also
been a lot of interest in leading contributions to other
important quantities, like, e.g., global fluctuations (i.e.,
variance of two traces) of the considered random
matrices. One should note that there are relations between
these two questions. We are not going to explore these
relations here, but we want to direct the reader's attention
to the so-called ``loop equations" in the physical
literature (see, e.g. \cite{Eyn}) and to the ``master
equation" in \cite{HT}.

We will concentrate in this paper on the second kind of
question.  As is well-known from the physical literature, in
many cases these leading contributions are given by planar
(or genus zero) diagrams and thus have quite a bit the
flavour of the combinatorics of free probability. In the
recent paper \cite{MN} this description was made precise for
the global fluctuations in the case of Wishart matrices, and
in particular the relevant set of planar diagrams (``annular
non-crossing permutations") was introduced and examined.
However, this description of the fluctuations in the large
$N$-limit was on a purely combinatorial level. Since it
looks quite similar to the description of free Poisson
distributions in terms of non-crossing partitions, one expects to
find some genuine free probability behind these
results. In particular, one would expect to have a description
on the level of operator algebras and to have
also a precise statement of the kind of ``freeness" that arises
here.

In this paper we will show that this is indeed the case. On
one hand, using the notion of a cyclic Fock
space, we can formulate the fluctuations in
terms of the usual creation, annihilation, and preservation
operators. On the other hand, we will also introduce an
abstract ``freeness" property for bilinear tracial
functionals, which not only give us a conceptual
understanding, but, on the other hand, is also crucial for
proving our main theorems on the fluctuations. 

Second order freeness, while stronger than the freeness of
Voiculescu, nevertheless appears to be a central feature of ensembles
of random matrices. Indeed, in this paper we prove that two standard
examples of random matrix ensembles exhibit second order freeness:
orthogonal families of Gaussian random matrices and orthogonal
families of Wishart random matrices are asymptotically free of
second order. Moreover in \cite{MSS} we show that independent Haar
distributed random unitary matrices are asymptotically free of
second order. 

The main results of the paper are thus. In section 5 (with proofs
in section 7) we show that semi-circular and compound Poisson
families on the full Fock space are free of second order. In
section 6 we establish the basic properties of second order
freeness and prove a general limit theorem.  In section 8 we
diagonalize the fluctuations in the Gaussian and Wishart case, thus
recovering and extending results of Cabanal-Duvillard \cite{C-D}.
In section 9 we prove asymptotic freeness of second order for
orthogonal families of Gaussian and Wishart random matrices.

\section{Preliminaries}

Here we collect some general notation and concepts which we
will use in the following.

Our presentation should be, by and large, self-contained,
however, it will rely of course on the basic ideas and
concepts of free probability.  For more details on this, one
should consult \cite{VDN, Voi2, NSp, HP}.  Furthermore, the
concepts of annular non-crossing  permutations and partitions
will play a crucial role. We will provide all relevant
information on them in the text. However, our presentation
will be quite condensed, and for further details one should
consult the original paper \cite{MN}.

\subsection{Some general notation}

For natural numbers $n,m\in\NN$ with $n<m$, we denote by
$[n,m]$ the interval of natural numbers between $n$ and $m$,
i.e.,
$$[n,m]:=\{n,n+1,n+2,\dots,m-1,m\}.$$
For a matrix $A=(a_{ij})_{i,j=1}^N$, we denote by $\Tr$ the
un-normalized trace and by $\tr$ the normalized trace,
$$\Tr(A):=\sum_{i=1}^N a_{ii},\qquad \tr(A):=\frac 1N \Tr(A).$$

For an $n\in\NN$, we will denote by $P(n)$ the set of
\emph{partitions} of $[1,n]$, i.e.,
$\sigma=\{B_1,\dots,B_r\}\in P(n)$ is a decomposition of
$[1,n]$ into disjoint subsets $B_i$: $B_i\not=\emptyset$ for
$i=1,\dots,r$, $B_i\cap B_j=\emptyset$ for $i\not= j$ and
$$[1,n]=\bigcup_{i=1}^r B_i.$$
The elements $B_i$ of $\sigma$ will be addressed as
\emph{blocks} of $\sigma$.

Given a mapping $i:[1,n]\to [1,N]$, the \emph{kernel},
$\ker(i)$, is defined as the partition of $[1,n]$, such that
two numbers $k,l\in[1,n]$ belong to the same block if and
only if $i(k)=i(l)$.

If we are considering classical random variables on some
probability space, then we denote by $\EE$ the expectation
with respect to the corresponding probability measure and by
$\kk_r$ the corresponding classical cumulants (as
multi-linear functionals in $r$ arguments); in particular,
$$\kk_1\{a\}=\EE\{a\}\qquad\text{and}\qquad
\kk_2\{a_1,a_2\}=\EE\{a_1a_2\}-\EE\{a_1\}\EE\{a_2\}.$$ 

\subsection{Annular non-crossing permutations and partitions}

The leading asymptotics of various random matrix quantities
can be described in terms of special ``planar" objects (see,
e.g., \cite{Eyn,Zvo}).  There are two equivalent ways of
formulating these results: a geometric ``genus"-expansion,
expressed by a sum over surfaces where the planar part
corresponds then to sums over surfaces of genus zero; an
algebraic description, where instead of using surfaces one
can sum over permutations and planarity is then a
geodesic-like condition on a length-function of these
permutations. If one prefers to associate partition like
pictures with permutations, then planarity is a condition
that these partitions have non-crossing diagrams (where,
however, one has to be careful about which drawings are
allowed).

We prefer to think in terms of permutations and
partitions. Let us recall the relevant definitions and
results.

Let, for $r\geq 1$, natural numbers $n(1),\dots,n(r)$ be
fixed. Consider a partition $\sigma\in
P(n(1)+\dots+n(r))$. In \cite{MN}, the class of
``multi-annular non-crossing partitions"
$\NC(n(1),\dots,n(r))$ was defined and, for $r=2$
(``annular" case), an extensive study of various
characterizations of such non-crossing partitions was made.
We will not go into details here, but we only want to state
the characterization which we will use. It will be the case
$r=2$ which is relevant for us; so let us use the notation
$\NC(n,m)$ for the $(n,m)$-annular non-crossing  partitions.
It is a good picture to think of two concentric circles,
with $n$ points on the outer and with $m$ points on the
inner. We put $[1,n]$ in clockwise order on the outer circle
and $[n+1,n+m]$ in counter-clockwise order on the inner
one. Adopting this convention will require that in some of
our formulas the indices corresponding to the outer circle
run in the opposite direction from the indices on the inner
circle.

\noindent\hfil\BoxedEPSF{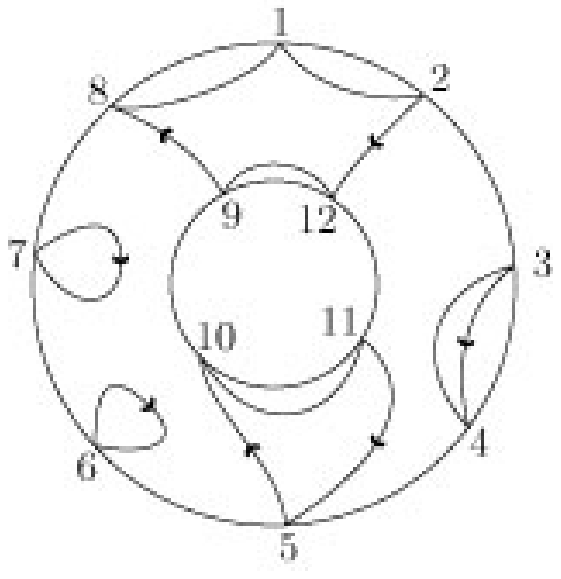}\hfil
\BoxedEPSF{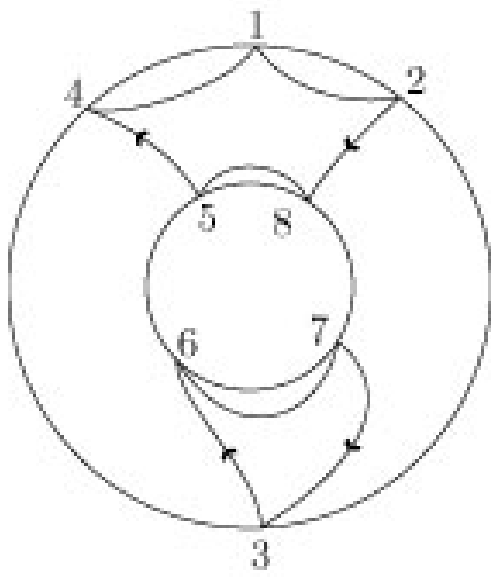}

\begin{quotation}{\small\noindent
{\bf Figure \figno.} On the left is the annular non-crossing permutation
$(1, 2, \ab 12, 9, 8)\, \ab (3, 4)\,\ab (5, 10, 11)\, (6)\, (7)$. On the
right is the permutation after the blocks that are contained in one
circle have been removed.
}\end{quotation}

Consider a $\sigma\in P(n+m)$.  We shall give a recursive
procedure for deciding if $\sigma$ is annular non-crossing.  
Suppose $\sigma$ has a block which is contained either in $[1,n]$ or
in $[n+1,n+m]$ and which consists of cyclically consecutive
numbers; then we remove this block and repeat the process
until we get a partition $\sigma' \in P(n'+m')$ with no
blocks 
which are contained in either $[1,n]$ or
$[n+1,n+m]$ and which consist of cyclically consecutive
elements.  Then,
by definition, $\sigma$ will be in $\NC(n,m)$ if and only if
$n',m'\geq 1$ and $\sigma'\in \NC(n',m')$. Thus it suffices
to say when $\sigma$ is in $\NC(n,m)$ for $\sigma$ with no
blocks which are contained in either $[1,n]$ or
$[n+1,n+m]$ and which consist of cyclically consecutive
elements. 

The characterizing property of such $\sigma$ is the
following: If we write the blocks $B\in\sigma$ in the form
$B=B'\cup B''$, where $B'\subset [1,n]$ and
$B''\subset[n+1,n+m]$, then, for all blocks $B$ of $\sigma$,
both parts, $B'$ and $B''$ must be non-empty and each of
them must consist of cyclically consecutive
numbers. Furthermore, the cyclic order of the restrictions,
$B_1', B_2', \dots , B_k'$, of the blocks of $\sigma$ to the
interval $[1, n]$ must be the reverse of the cyclic order of
the restrictions, $B_1'', B_2'', \dots B_k''$, of the blocks
to the interval $[n+1,n+m]$.  Note that this
characterization contains the statement that a
$\sigma\in\NC(n,m)$ is connected in the sense that at least
one block of $\sigma$ contains elements both from $[1,n]$
and from $[n+1,n+m]$ (i.e., $\sigma$ connects the two
circles).

In the context of random matrices, it is permutations, not
partitions, which appear in calculations.  In order to go
over from partitions $\sigma$ to permutations $\pi$ one has
to choose a cyclic order on each of the blocks of $\sigma$.
Choosing such an order for each block will produce an
``annular non-crossing permutation" out of an annular
non-crossing partition. The set of annular non-crossing
permutations is denoted by $\SNC(n,m)$ -- and by
$\SNC(n(1),\dots,n(r))$ in the multi-annular case -- and it
is this set which was the main object of interest in
\cite{MN}. In [MN, \S 6] it was shown that a
permutation $\tau$ is in $S_{NC}(n,m)$ if two conditions are
satisfied. The first condition is {\it connectedness}: at
least one cycle of $\tau$ connects the two circles and the
second is {\it planarity}: the {\it geodesic condition} must
be satisfied: $\#(\pi) + \#(\pi^{-1} \gamma) = m + n$, where
$\#(\pi)$ denotes the number of cycles of $\pi$ and $\gamma$
is the permutation with two cycles: $\gamma = \gamma_{n,m} =
(1,2 , 3 \dots , n)(n+1, \dots, n+m)$

We wish to describe what it means to be a a non-crossing
permutation on an $r$ multi-annulus. By an $r$
multi-annulus we mean a collection of $r$ circles with
$n(1)$ points on the first circle, $n(2)$ points on the
second circle, \dots , $n(r)$ points on the $r^{th}$
circle. Connectivity of $\tau$ means that every pair of
circles is connected by at least one cycle of $\tau$. The
planarity of $\tau$ is defined using a geodesic condition.
Let $\gamma_{n(1), n(2), \dots n(r)}$ be the permutation of
$[n(1) + \cdots + n(r)]$ with $r$ cycles --- the $r^{th}$
cycle being $(n(1) + \cdots n(r-1) + 1, \dots , n(1) +
\cdots n(r))$. $\tau$ will be planar if $\tau$ satisfies
the geodesic condition $\#(\tau) + \#(\gamma_{n(1), n(2),
\dots n(r)} \tau^{-1}) = n(1) + \cdots + n(r) + 2 - r$.

As observed in \cite{MN}, there is not necessarily a unique
choice of a cyclic order on a block of $\sigma$; to put it
another way, the mapping from $\pi$ to $\sigma$ (which
consists in forgetting the order on the cycles) is not
injective.  However, this deviation from injectivity is not
too bad.  Let us consider a block $B\in\sigma$, and denote
by $B':=B\cap [1,n]$ and $B'':=B\cap [n+1,n+m]$ the parts of
$B$ lying on the first and on the second circle,
respectively. On each of the two circles we respect the
given cyclic order on $(1,\dots,n)$ and on $(n+1,\dots,n+m)$
The allowed orders on $B$ thus consist of choosing a `first'
element of $B'$ and a `first' element of $B''$; then the
order on $B$ is obtained by running through $B'$ from the
first to the last element, then going over to the first
element in $B''$ and continuing in $B''$ to the last
element.  Hence the only freedom we have is the choice of
first elements in $B'$ and in $B''$.

Let us call a block $B\in\sigma$ a \emph{through-block}, if
both $B':=B\cap [1,n]$ and $B'':=B\cap [n+1,n+m]$ are
non-empty. Then only $\sigma$ with exactly one through-block
have two or more $\pi$'s in $\SNC(n,m)$ as preimages.
Namely, if $B=B'\cup B''$ is the unique through-block of
such a $\sigma$, then every element from $B'$ can be chosen
as first element, and the same for $B''$, thus there are
exactly $\vert B'\vert\cdot \vert B''\vert$ possible choices
of cyclic orders for $B$.  If, however, there is more than
one through-block, then the first element on each component
of them is uniquely determined and there is exactly one
possible order for each block.

For Gaussian random matrices only non-crossing pairings will
play a role.  These are those annular non-crossing
partitions for which each block consists of exactly two
elements.  One should note that in this case the distinction
between permutations and partitions vanishes, because for
pairings there is always exactly one possibility for putting
an order on blocks. We will denote the set of annular
non-crossing pairings by $\NC_2(n,m)$; and, for the
multi-annular situation, by $\NC_2(n(1),\dots,n(r))$. 
In the multi-annular case the geodesic condition can be written
$\#(\gamma_{n_1, \dots, n_r} \pi) = 2 - r + (n_1 + \cdots +
n_r)/2$.

\section{Combinatorial description of global fluctuations}

We are interested in the fluctuations of various types of
$N\times N$ random matrices around their large $N$-limit.
Here, we are going to consider two classes of random
matrices, namely Gaussian random matrices and (a
generalization of) Wishart matrices. Let us fix the notation
for our investigations.

\subsection{Semi-circular case}\label{gaussianrm}

Let $(X_N)_{N\in\NN}$ be a sequence of $N\times N$ Hermitian
Gaussian random matrices. Then, in the limit $N\to\infty$,
$X_N$ converges to a semi-circular variable $s$.  Let us
consider directly the case of several such Gaussian random
matrices. The entries of different random matrices need not
be independent from each other, but they have to form a
Gaussian family. A convenient way to describe such a
situation is to index the matrices by elements from some
real Hilbert space $\HH_\RR$, such that the covariance
between entries from $X_N(f)$ and $X_N(g)$ is given by the
inner product $\la f,g\ra$.  More precisely we say that
$\{X_n(f) \}_{f \in \HH_\RR}$ is a {\it family} of Hermitian
Gaussian random matrices if $X_N(f) = (x_{i,j}(f))_{i,j = 1}^N$
and the entries $\{ x_{i,j}(f) \mid 1 \leq i, j \leq N, f \in
\HH_\RR \}$ form a Gaussian family with covariance given by 
$$\EE\{x_{ij}(f)x_{kl}(g)\}=0
\mbox{ for } i < j, k < l, \mbox{ and } f,g \in \HH_\RR \mbox{ and }$$
$$ \EE\{x_{ij}(f)\bar
x_{kl}(g)\}=\delta_{ik}\delta_{jl}\cdot\frac 1N \la
f,g\ra \mbox{ for } i \leq j, k \leq l, \mbox{ and } f,g \in \HH_\RR$$

By Wick's formula (see e.g. \cite[Sections 1.3 and 1.4]{J}) we have 
\[
E\{x_{i_1, j_1}(f_1) x_{i_2, j_2}(f_2) \cdots x_{i_{2k},
j_{2k}}(f_{2k})\}  =
\sum_{\pi \in \cP_2(2k)} \prod_{(r,s) \in \pi} 
\langle
x_{i_r, j_r}(f_r), x_{i_s, j_s}(f_s) \rangle
\]
where the sum is over all pairings $\pi$ of $[2k]$ and the
contribution of each pairing is the product of $\langle
x_{i_r, j_r}, x_{i_s, j_s} \rangle$ over all pairs $(r, s)$ of
$\pi$.

Then, Voiculescu's extension of
Wigner's theorem to this multi-dimensional case states that,
for $N\to\infty$, such a family of random matrices converges
to a semi-circular system with the same covariance. We want
to look more closely on that convergence and investigate the
``global fluctuations" around this semi-circular limit; this
means, we want to understand the asymptotic behaviour of
traces of products of our random matrices.  It turns out
that, with the right scaling with $N$, these random
variables converge to a Gaussian family and thus the main
information about them is contained in their covariance.

If one invokes the usual genus expansions for expectations
of Gaussian random matrices then one gets quite easily the
following theorem. It turns out that the leading orders are
given by planar pairings. Since we are looking on cumulants
and not just moments, the relevant pairings also have to
connect the $r$ circles.

\begin{theorem}\label{hermitianmatrix}

Let $X_N(f)$ ($f\in\HH_\RR$) be a family of Hermitian
Gaussian random matrices. Let $\kk_r$ denote the $r^{th}$
classical cumulant (considered as multi-linear mapping of
$r$ arguments), then for $f_1,\dots,f_{n_1+\cdots+
  n_r}\ab  \in \ab \HH_\RR$, the leading order of the cumulants of
the random variables
\begin{multline}
\big\{\Tr\big(X_N(f_1)\cdots X_N(f_{n(1)})\big),\ab
\dots ,\\
\Tr\big( X_N(f_{n(1) + \cdots + n(r-1) + 1}) \cdots
X_N(f_{n(1) + \cdots + n(r)})\big)\big\}
\end{multline}
are given by
\begin{multline}\label{gaussian-cumulants}
\kk_r\Bigl\{\Tr[X_N(f_1)\cdots X_N(f_{n(1)})],\cdots, \\
\Tr[X_N(f_{n(1)+\dots+n(r-1)+1})\cdots X_N(f_{n(1)+\dots
+n(r)})]\Bigr\}\\ =N^{2-r}\cdot \sum_{\pi\in NC_2(n(1),\dots,n(r))}
\prod_{(i,j)\in\pi}\la f_i,f_j\ra + O(N^{-r}).
\end{multline}
\end{theorem}

\noindent{\it Proof:}
Let 
$n = n(1) + \cdots + n(r)$, $\gamma$ be the permutation of $[n_1 + \cdots
+ n_r]$ with the $r$ cycles $(1, \dots, n_1) (n_1 + 1, \dots , n_1 + n_2)
\cdots (n_1 + \cdots + n_{r-1}+1, \cdots + n_1 + \cdots + n_r)$, and
$Y_i =
\tr\big(X_N(f_{n(1) + \cdots + n(i-1)+1}) \ab\cdots
X_N(f_{n(1) + \cdots + n(i)})\big)$. Then by Wick's formula

\begin{eqnarray*}
\EE(Y_1 \cdots  Y_r)
&=&
\sum_{i_1, \dots , i_n = 1}^N
\EE( x_{i_1, i_{\gamma(1)}}(f_1) \cdots x_{i_n, i_1}(f_n))
\\
&=&
\sum_{i_1, \dots , i_n = 1}^N
\sum_{\pi \in \cP_2(n)} \prod_{(k,l) \in \pi}
\langle f_k, f_l \rangle
\delta_{i_k, i_{\gamma(l)}} \delta_{i_l, i_{\gamma(k)}} \\
&=&
\sum_{\pi \in \cP_2(n)} \prod_{(k,l) \in \pi}
\langle f_k, f_l \rangle
\sum_{i_1, \dots , i_n = 1}^N
\delta_{i_k, i_{\gamma(l)}} \delta_{i_l, i_{\gamma(k)}} \\
&=&
\sum_{\pi \in \cP_2(n)} \prod_{(k,l) \in \pi}
\langle f_k, f_l \rangle N^{\#(\gamma \pi)}
\end{eqnarray*}

Following the argument in [MN, proof of Proposition
9.3] we have

\[
k_r(Y_1, \dots , Y_r) =
\mathop{\sum_{\pi \in \cP_2(n)} }_{%
\pi \textrm{ is connected}}
\prod_{(k,l) \in \pi}
\langle f_k, f_l \rangle
N^{\#(\gamma \pi)}
\]

The terms of highest order are the planar ones thus we
obtain equation (\ref{gaussian-cumulants}).\qed

\medskip

This theorem contains all relevant combinatorial information
about the asymptotic behaviour of our traces. Since an
increase of the number of arguments of the cumulants
corresponds to a decrease in the exponent of $N$, a cumulant
$\kk_r$ will always dominate a cumulant $\kk_p$ if $r<p$. So
in leading order only the first cumulant survives in the
limit, which gives us the following statement analogous to
the law of large numbers.

\begin{corollary}

For each $f_1,\dots,f_n\in\HH_\RR$, the random variables
$$\Big\{\tr\Big(\,X_N(f_1)\cdots X_N(f_n)\,\Big)\Big\}_N$$
converge in distribution to the constant random variables
$\alpha(f_1,\dots,f_n)\cdot 1$, where
$$\alpha(f_1,\dots,f_n)=\sum_{\pi\in
NC_2(n)}\prod_{(i,j)\in\pi} \la f_i,f_j\ra.$$

\end{corollary}

This corollary is of course just a reformulation of
Voiculescu's result that $\EE\{\tr[X_N(f_1)\cdots
  X_N(f_n)]\}$ has, in the limit $N\to\infty$, to agree with
the corresponding moments of a semi-circular family.

But we can now go a step further. If we subtract the mean of
the random variables, then the first cumulants are shifted
to zero and it will be the second cumulants which survive --
after the right rescaling. Since higher cumulants vanish
compared to the second ones, we get Gaussian variables in
the limit.

\begin{corollary}\label{gaussiansecondorder}

Consider the (magnified) fluctuations around the limit
expectation,
\begin{align*}
  F_N(f_1,\dots,f_n)&:=N\cdot\bigl(\tr[X_N(f_1)\cdots
  X_N(f_n)]- \alpha(f_1,\dots,f_n)\bigr)\\
& \mbox{} =\Tr[X_N(f_1)\cdots X_N(f_n)]-N\alpha(f_1,\dots,f_n).\notag
\end{align*}
The family of all fluctuations
$\bigl(F_N(f_1,\dots,f_n)\bigr)_{n\in \NN,f_i\in\HH_\RR}$
converges in distribution towards
$\bigl(F(f_1,\dots,f_n)\bigr)_{n\in \NN,f_i\in\HH_\RR}$, a
centered Gaussian family with covariance given by
\begin{equation*}
\EE\{F(f_1,\dots,f_n)\cdot F(f_{n+1},\dots,f_{n+m})\}=
\sum_{\pi\in NC_2(n,m)}\prod_{(i,j)\in\pi}\la f_i,f_j\ra.
\end{equation*}

\end{corollary}

Our goal now is to present a conceptual understanding of
this form of the covariance; in particular one that would
easily diagonalize it. In principle, this is a purely
combinatorial problem. However, our point of view is that
limits of random matrices which have the flavour of free
combinatorics should also have a description in terms of the
operator side of free probability, i.e., operators on full
Fock spaces. We will provide such a description and show
that it diagonalizes our covariance.

\subsection{Compound Poisson case}
Let $(X_N)_{N\in\NN}$ be a sequence of $N\times N$ complex
Gaussian random matrices (i.e.  the entries of $X_N$ are
independent centered complex Gaussians with variance $1/N$)
and let $(D_N)_{n\in\NN}$ be a sequence of $N \times N$ non-random
matrices for which a
limit distribution exists as
$N\to\infty$. Then, in the limit $N\to\infty$,
$\{X_N,X_N^\ast, D_N\}$ converges in distribution to
$\{c,c^\ast, d\}$, where $c$ is a circular element, $d$ has
the limit distribution of the $D_N$, and $\{c, c^\ast \}$
and $d$ are free. In particular, $X_N^\ast D_N X_N$
converges to $c^\ast dc$, which is a free compound Poisson
element, see \cite[4.4]{Sp2}. We shall discuss the
fluctuations of the random matrices
$$P_N:=X_N^\ast D_N X_N$$
around the limit $c^\ast dc$.  Since $P_N$ is a generalization of
a Wishart matrix, we will call it in the following a
\emph{compound Wishart matrix}.

As we shall see it is appropriate to consider a more general
situation.  Namely, consider not just a single non-random
matrix $D_N$, but also all its powers $D_N^k$ at the same
time, or more generally, let us consider a family $\{
D_1^{(N)}, D_2^{(N)}, D_3^{(N)}, \dots, D_p^{(N)}\}_N$ of $N
\times N$ complex matrices. We shall say the family {\it
converges in distribution} if there are operators $d_1, d_2,
d_3,\dots , d_p$ and a tracial state $\psi$ on $\cC$, the
complex *-algebra generated by $\{1, d_1, d_2, \dots ,
d_p\}$, such that
$$\lim_{N\to\infty}\tr[D_{i_1}^{(N)}\cdots
D_{i_k}^{(N)}]=\psi(d_{i_1} \cdots d_{i_k})$$
for all  $i_1, i_2, \dots , i_k$.

We are again interested in global fluctuations of these
matrices in the limit $N\to\infty$; i.e., we want to
consider the asymptotic behaviour of mixed moments of our
random matrices. Again the key point is the understanding
of the leading order of the cumulants in these traces. 
This leading order is described by summing over
non-crossing permutations, but in contrast to the
semi-circular case, all permutations contribute, not just
pairings. In order to describe the contribution of such a
general non-crossing permutation, we need the following
notation.

\begin{notation}
Let $(\cA, \psi)$ be a unital algebra with a tracial state
$\psi$; for each $\pi \in S_p$ we shall define a $p$-linear
functional, $\psi_\pi$, on $\cA \times \cdots \times \cA$. Write
$\pi = c_1 \cdot c_2 \cdot \cdots \cdot c_k$ as a product of
disjoint cycles, and for each $i$, $c_i = (r_{i,1}, \dots ,
r_{i, l_i})$.  Then define the $p$-linear
functional $\psi_\pi$ by
\[
   \psi_\pi(a_1, a_2, a_3, \dots , a_p) = \prod_{i=1}^k
   \psi(a_{r_{i,1}} \cdots a_{r_{i,l_i}})
\]
\end{notation}

Note that we need $\psi$ to be a trace, because a cycle $c$
comes only with a cyclic order.

An example of this notation is the following, take
$$\pi=\{(1,2,6),(3,4,5)\}\in \SNC(3,3).$$
Then
$$\psi_\pi(a_1,a_2,a_3,a_4,a_5,a_6)=\psi(a_1a_2a_6)\cdot
\psi(a_3a_4a_5).$$
Note also that the cyclic order is important. In $\SNC(2, 1)$
consider
$$\pi_1=\{(1,2,3)\}\qquad\text{and}\qquad\pi_2=\{(1,3,2)\}.$$
Although their block structure is the same, as permutations
they are different elements from $\SNC(2, 1)$ and we have
$$\psi_{\pi_1}(a_1,a_2,a_3)=\psi(a_1a_2a_3)
\qquad\text{and}\qquad
\psi_{\pi_2}(a_1,a_2,a_3)=\psi(a_1a_3a_2).
$$
We shall denote the number of cycles in the permutation $\pi$ by
$\#(\pi)$.

Let us now state the basic combinatorial description of the
leading order of cumulants in traces of products of our
compound Wishart matrices. For the usual Wishart matrices
this was derived in \cite{MN}.  Our more general version
follows by the same kind of calculations ({\it c.f.}
Capitaine and Casalis \cite[\S 5]{CC}).

\begin{theorem}

Let $\{ X_N\}_N$ be a sequence of complex
Gaussian random matrices. Put
$$P_N(D_i):=X_N^\ast D_i^{(N)}X_N$$
Let $(\kk_r)_{r\in\NN}$ denote the classical cumulants, then
we have for all $r\in\NN$
\begin{multline} \label{wishartcumulant}
\lim_N N^{r-2}
\kk_r\Bigl\{\Tr[P_N(D_1)\cdots P_N(D_{n_1})],\cdots, \\
\Tr[P_N(D_{n_1 + \cdots + n_{r-1}+1}) \cdots  
P_N(D_{n_1 + \cdots + n_r})]\Bigr\}\\ 
=
\sum_{\pi\in \SNC(n(1),\dots,n(r))}
\psi_\pi(d_1, \dots,d_{n_1 + \cdots + n_r}) \\
\end{multline}

\end{theorem}

\noindent{\it
Proof:}
Let $n = n_1 + \cdots + n_r$ and let $\gamma$ be the permutation with
$r$ cycles: $(1, \dots, n_1) (n_1 + 1, \dots, n_1 + n_2) \cdots (n_1
+ \cdots + n_{r-1} +1, \dots, n_1 + \cdots + n_r)$. Let
\[Y_l = \Tr(P_N(D_{n_1 + \cdots + n_{l-1}+1}) \cdots\ab  
P_N(D_{n_1 + \cdots + n_l}))\] By \cite[Theorem 2]{GLM}
\begin{equation}\label{fullexpansion}
\EE(Y_1 \cdots Y_r)  =
\sum_{\sigma \in S_n} N^{\#(\sigma^{-1}\gamma)-n}
\Tr_{\sigma}(D_1, D_2, \dots , D_n)
\end{equation}

For $\sigma \in S_n$ let $\sigma \vee
\gamma$ be the partition of $[n]$ whose blocks are the orbits of
the group generated by $\sigma$ and $\gamma$. $\sigma \vee
\gamma$ also defines a partition of the $r$ cycles of $\gamma$.
Let us denote this partition of $[r]$ by $A_\gamma(\sigma)$.
Conversely let $s$ be the number of cycles of $\sigma$, $\sigma
\vee \gamma$ determines a partition of the cycles of $\sigma$;
we shall denote this by $A_\sigma(\gamma)$. Note that if
$A_\gamma(\sigma) = 1_r$ then $\sigma \vee \gamma = 1_n$ and
thus $A_\sigma(\gamma) = 1_s$.

For a partition $A = \{A_1, \dots , A_k\}$ of $[r]$ let
\[\EE_A(Y_1, \dots , Y_r) = \prod_{l=1}^k \EE( \prod_{i \in A_l}
Y_i)\] If $A = 1_r$ then $\EE_A(Y_1, \dots , Y_r) = \EE(Y_1
\cdots Y_r)$. Equation (\ref{fullexpansion}) can now be extended
easily to obtain
\[
\EE_A(Y_1, \dots , Y_r) =
\mathop{\sum_{\sigma \in S_n}}_{A_\gamma(\sigma) \leq A}
N^{\#(\sigma^{-1}\gamma)-n}
\Tr_\sigma(D_1, \dots , D_n)
\]

Let $\mu(A, B)$ be the M\"obius function of the lattice
of partitions; in particular $\mu(A, 1_r) = (-1)^{\#(A)-1}
(\#(A)-1)!$. Note that 
\[
\sum_{B \leq A \leq 1_r}
\mu(A, 1_r) =
\begin{cases}
1 & B = 1_r \\
0 & B < 1_r \\
\end{cases}
\]

\begin{eqnarray*}\lefteqn{
k_r(Y_1, \dots , Y_r) }\\
&=&
\sum_{A \in \cP(r)} \mu(A, 1_r)
\EE_A(Y_1, \dots , Y_r) \\
&=&
\sum_{A \in \cP(r)} \mu(A, 1_r)
\mathop{\sum_{\sigma \in S_n}}_{A_\gamma(\sigma) \leq A} 
N^{\#(\sigma^{-1}\gamma)-n}
\Tr_\sigma(D_1, \dots , D_n) \\
&=&
\sum_{\sigma \in S_n}
N^{\#(\sigma^{-1}\gamma)-n}
\Tr_\sigma(D_1, \dots , D_n)  
\mathop{\sum_{A \in \cP(r)}}_{A_\gamma(\sigma) \leq A}  \mu(A, 1_r) \\
&=&
\mathop{\sum_{\sigma \in S_n}}_{\sigma \vee \gamma = 1_n}
N^{\#(\sigma^{-1}\gamma)-n}
\Tr_\sigma(D_1, \dots , D_n) \\
&=&
\mathop{\sum_{\sigma \in S_n}}_{\sigma \vee \gamma = 1_n}
N^{\#(\sigma^{-1}\gamma)-n}
\Tr_\sigma(D_1, \dots , D_n)
\end{eqnarray*}

Recall that for $\sigma \in S_n$ with $\sigma \vee \gamma =
1_n$ there is an integer $g = g(\sigma)$ such that $\#(\sigma) +
\#(\sigma^{-1} \gamma) + \#(\gamma) = n + 2(1 - g)$
and that $\sigma \in S_{NC}(n_1, \dots , n_r)$ means
that $\sigma \vee \gamma = 1_n$ and $g(\sigma) = 0$.

\begin{eqnarray*}\lefteqn{
k_r(Y_1, \cdots , Y_r) } \\
&=& N^{-n}
\mathop{\sum_{\sigma \in S_n}}_{\sigma \vee \gamma = 1_n}
N^{\#(\sigma^{-1}\gamma)}\, N^{\#(\sigma) }\,
\tr_\sigma(D_1, \dots , D_n) \\
&=&
\mathop{\sum_{\sigma \in S_n}}_{\sigma \vee \gamma = 1_n}
 N^{2 - r - 2g(\sigma)}
\tr_\sigma(D_1, \dots , D_n)  \\
&=&
N^{2-r} \kern -1.5em
\sum_{\sigma \in S_{NC}(n_1, \dots , n_r)} \kern -1em
\tr_\sigma(D_1, \dots , D_n)  
+ O(N^{-r})
\end{eqnarray*}

Since $\lim_N \tr_\sigma(D_1, \dots , D_n)  =
\psi_\sigma(d_1, \dots, d_n)$ we have the required result. \qed

\medskip

This theorem contains again all relevant information about
the limit behaviour of the random variables
$\Tr(P_N(D_{i_1})\cdots P_N(D_{i_n}))$.  First, we have the
following statement analogous to the law of large numbers.

\begin{corollary}

The random variables $\{\tr[P_N(D_{i_1})\cdots
P_N(D_{i_n})]\}_{i_1, \dots i_n}$ converge in distribution
to  constant random variables $\beta(d_{i_1},\ab \dots,\ab
d_{i_n})\cdot 1$, where
\begin{equation}\label{mom:poisson}
\beta(d_{i_1},\dots,d_{i_n}):= 
\sum_{\pi\in NC(n)}\psi_\pi(d_{i_1},\dots,d_{i_n}).
\end{equation}

\end{corollary}

The form of $\beta(d_{i_1},\dots,d_{i_n})$ is, of course, in
agreement with the fact that $\EE\{\tr[P_N(D_{i_1})\cdots
P_N(D_{i_n})]\}$ has, in the limit $N\to\infty$, to agree
with the corresponding moment of the compound free Poisson
variables,
\[\psi(c^\ast d_{i_1} c \cdot c^\ast d_{i_2} c \cdots
c^\ast d_{i_n} c)\] where $c$ is a circular random variable
*-free from $\{d_1, \dots, d_p\}$.  Again, we magnify the
fluctuations around that limit, thus shifting the first
cumulants to zero and getting only a non-vanishing limit for the
second cumulants -- hence getting normal limit fluctuations.

\begin{corollary}\label{wishartsecondorder}

Consider the (magnified) fluctuations around the limit
value,
\begin{align}
F_N(D_{i_1},\dots,D_{i_n}):&=N \cdot \Big(\tr[P_N(D_{i_1})\cdots
P_N(D_{i_n})]- \beta(d_{i_1},\dots,d_{i_n})\Big)\\ 
&=\Tr[P_N(D_{i_1})\cdots P_N(D_{i_n})]-N\beta(d_{i_1},\dots,d_{i_n}).\notag
\end{align}

The family of all fluctuations $(F_N(D_{i_1},\dots,D_{i_n} )
)_{n \in \NN}$ converges in distribution towards a centered
Gaussian family $\Big( (F(d_{i_1},\dots,d_{i_n})\Big)_{i_1,
\dots , i_n}$, with covariance given by
\begin{multline}
\EE\{F(d_{i_1},\dots,d_{i_m})\cdot F(d_{i_{m+1}},\dots,d_{i_{m+n}})\} \\
=\sum_{\pi\in\SNC(m,n)}\psi_\pi(d_{i_1}, \dots, d_{i_{n+m}}).
\end{multline}
\end{corollary}

Again, it remains to understand this covariance and we will
be aiming at a more operator-algebraic description of these
fluctuations in order to attack this combinatorial problem.

\section{Realization of semi-circular and free compound Poisson 
elements on Fock spaces}

The main theme of our investigations is the conviction that
wherever planar or non-crossing objects arise, there is some
free probability lurking behind the picture. Since the
fluctuations of our Gaussian and Wishart random matrices can
be described combinatorially in terms of non-crossing
permutations, we expect also some operator-algebraic or some
more abstract ``free" description of this situation. Our
main results in the coming sections will provide these
descriptions. Let us begin by recalling the realization of a
semi-circle and a compound Poisson distribution on a full
Fock space by using creation, annihilation, and preservation
operators.

\subsection{Semi-circular case}

For a real Hilbert space $\HH_\RR$ with complexification
$\HH$, we consider the full Fock space
$$\cF(\HH):=\bigoplus_{n=0}^\infty\HH^{\otimes n}=
\CC\Omega\oplus \HH\oplus \HH^{\otimes 2} \oplus\dots$$
and define, for $f\in\HH$, the \emph{creation operator}
$l(f)$ by
$$l(f)\Omega=f$$ and
$$l(f)f_1\odo f_n=f\otimes f_1\odo f_n$$
and the 
\emph{annihilation operator} $l^*(f)$,
by 
$$l^*(f)\Omega=0$$
and
\begin{align*}
l^*(f)f_1\odo f_n&=\la f_1,f\ra f_2\odo f_n\\
\end{align*} 
($n\in\NN$, $f_1,\dots,f_n\in\HH$).

For $f\in\HH_\RR$, we put
$$\omega(f):=l(f)+l^*(f)$$
and we will denote by $\cA(\HH_\RR)$ the complex unital
$*$-algebra generated by all $\omega(f)$ for
$f\in\HH_\RR$. Note that all $\omega(f)$ are self-adjoint
and that the vector $\Omega$ is cyclic and separating for
the algebra $\cA(\HH_\RR)$.

If we define on $\HH$ an involution $f\mapsto \bar f$ by
$$\overline{f_1+if_2}:=f_1-if_2\qquad
\text{for $f_1,f_2\in\HH_\RR$},$$
then $f\mapsto \omega(f)$ extends from a real linear mapping
on $\HH_\RR$ to a complex linear mapping on $\HH$ with
$$\omega(f)=l(f)+l^*(\bar f)\qquad (f\in\HH).$$
Note that the unital $*$-algebra generated by all
$\omega(f)$ with $f\in\HH$ is just $\cA(\HH_\RR)$.

It is well known (see, e.g., \cite{VDN}) that these
operators $\omega(f)$ have a semi-circular distribution and
thus the asymptotics of the expectation of traces of
Gaussian random matrices can also be stated as follows.

\begin{proposition}

Let $X_N(f)$ ($f\in\HH_\RR$) be a family of Hermitian
Gaussian random matrices.  Then for all
$f_1,\dots,f_n\in\HH_\RR$
\begin{equation}
\lim_{N\to\infty} \EE\{\tr[X_N(f_1)\cdots X_N(f_n)]\}
=\la\omega(f_1)\cdots \omega(f_n)\Omega,\Omega\ra.
\end{equation}

\end{proposition}

Let us in this context also recall the definition of the
Wick products.

\begin{definition}
For $f_1,\dots,f_n\in\HH$ the \emph{Wick product} $W(f_1\odo
f_n)$ is the unique element of $\cA(\HH_\RR)$ such
\begin{equation}
W(f_1\odo f_n)\Omega=f_1\otimes\dots\otimes f_n.
\end{equation}
For $n=0$, this has to be understood as $W(\Omega)=1$.
\end{definition}

Since $\Omega$ is cyclic and separating for $\cA(\HH_\RR)$,
these Wick products exist and are uniquely determined.

From the definition of the creation and annihilation
operators it is clear that these Wick products satisfy for
all $f,f_1,\dots,f_n\in\HH$ the relation
$$\omega(f)W(f_1\odo f_n)=W(f\otimes f_1\odo f_n)+\la f_1,
\bar f\ra W(f_2\odo f_n).$$
This can also be used as a recursive definition for the Wick
products and shows that $W(f_1\odo f_n)$ is a polynomial in
$\omega(f_1),\dots,\omega(f_n)$.

In the case $f=f_1=\dots=f_n$ this reduces to the three-term
recurrence relation for the Chebyshev polynomials and shows
that
$$W(f^{\otimes n})=U_n(\omega(f)/2),$$
where $U_n$ is the $n^{th}$ Chebyshev polynomial of the
second kind.

\subsection{Compound Poisson case}

In this case we start with a unital $*$-algebra $\cC$
equipped with a tracial state $\psi$ and represent $\cC$,
via the GNS-representation, on $\HH:=
\overline{\cC}^{\la\cdot,\cdot\ra}$, where the inner product
on $\HH$ is given by
$$\la d_1,d_2\ra:=\psi(d_2^*d_1).$$

Then we take the full Fock space $\cF(\HH)$ and consider
there as before the creation and annihilation operators
$l(d)$ and $l^*(d)$, respectively.  But now we have, in
addition, also to consider the \emph{preservation} (or
\emph{gauge}) \emph{operator} $\Lambda(d)$ ($d\in\cC$) which
is defined by
$$\Lambda(d)\Omega=0$$
and 
$$\Lambda(d) f_1\otimes\cdots\otimes f_n= (df_1)\otimes f_2\otimes
\cdots\otimes f_n$$
for $f_1,\dots,f_n\in\HH$. Note that the multiplication
$\cC\times \cC\to\cC$ extends to a module action
\begin{align*}
\cC\times \HH&\to \HH\\
(d,f)&\mapsto df.
\end{align*}

For $d\in\cC$ we define now
\begin{equation}
p(d):=l(d)+l^*(d^*)+\Lambda(d)+\psi(d)1,
\end{equation}
and we will denote by $\cA(\cC)$ the unital $*$-algebra
generated by all $p(d)$ for $d\in\cC$. Note that we have
$$p(d)^*=p(d^*)\qquad\text{for $d\in\cC$}.$$
 
One knows (see, e.g., \cite{GSS,Sp1,NSp}) that these
operators $p(d)$ give a realization of compound Poisson
elements, i.e., their moments are given by
Eq. (\ref{mom:poisson}).  Thus we can state the asymptotics
of the expected value of traces of our compound Wishart
matrices also in the following form.

\begin{proposition}

Suppose the family $\{D^{(N)}_1, \dots , D_p^{(N)}\}$
converges in distribution to $\{d_1, \dots , d_p\}$ in
$(\cC,\psi)$.  Let $(X_N)_{N\in\NN}$ be a sequence of
$N\times N$ Hermitian Gaussian random matrices. Then
\begin{equation}
\lim_{N\to\infty} \EE\{\tr[P_N(D_{i_1})\cdots
P_N(D_{i_n})]\} =\la p(d_{i_1})\cdots
p(d_{i_n})\Omega,\Omega\ra.
\end{equation}
\end{proposition}

Again, Wick products will play a role in this context. As
before, these should be polynomials in the $\{p(d)\mid
d\in\cC\}$ with the defining property that
$$W(d_1\otimes\dots\otimes d_n)\Omega=d_1\otimes\dots\otimes
d_n.$$ 
However, in contrast to the semi-circular case, the
multiplication of $d$'s in the arguments (under the action
of $\Lambda$) has the effect that in order to produce
counter-terms for $p(d_1)\cdots p(d_n)\Omega$ to get
$d_1\otimes \dots\otimes d_n$ one also has to involve
operators like $p(d_1d_2)$ etc. This means that
$W(d^{\otimes n})$ is in general not just a polynomial in
$p(d)$, but some polynomial in all $\{p(d^k)\mid k\leq
n\}$. In particular, in general there is no relation between
Wick polynomials $W(d^{\otimes n})$ and the orthogonal
polynomials with respect to the distribution of $p(d)$.
From the point of view of Levy processes the occurrence of
$p(d^k)$ in $W(d^{\otimes n})$ is not very surprising,
because this corresponds to the higher diagonal measures
(Ito-formulas) and a Levy process should come along with its
higher variations.

It appears that Anshelevich \cite{Ans} was the first to
introduce and investigate these polynomials in this
generality (and also some $q$-deform\-ations thereof). Since
these polynomials appear implicitly in the classical case in
a paper of Kailath and Segall, he called them free
Kailath-Segall polynomials.

By taking into account the action of our operators on the
full Fock space one sees quite easily that these Wick
products should be defined as follows.

\begin{definition}
For a given algebra $\cC$ with state $\psi$, the \emph{Wick
  products} or \emph{free Kailath-Segall polynomials} of the
corresponding compound Poisson distribution are recursively
defined by ($d,d_1,\dots,d_n\in\cC$)
$$W(d)=p(d)-\psi(d)1$$
and
\begin{align*}
W(d\otimes d_1\otimes\cdots \otimes d_n)&=
p(d) W(d_1\otimes d_2\otimes\dots\otimes d_n)\\
&\qquad-
\psi(dd_1)W(d_2\otimes d_3\otimes\dots\otimes d_n)\\
&\qquad-
W(dd_1\otimes d_2\otimes\cdots\otimes d_n)\\
&\qquad-\psi(d) W(d_1\otimes d_2\otimes\cdots\otimes d_n)
\end{align*}    
\end{definition}

\section{Cyclic Fock space}

Our main aim is to express the formulas for the limit
fluctuations of Gaussian random matrices and of compound
Wishart matrices also with the help of the operators
$\omega(f)$ and $p(d)$, respectively.  In order to do so we
have, however, to introduce another variant of a Fock
space. Whereas the elements in the full Fock space,
$f_1\otimes\dots\otimes f_n$, are linear kind of objects --
with a beginning and an end - we are looking on traces and
thus should identify the beginning and the end in a cyclic
way.

Here are two versions of such a cyclic Fock space, the first
one over arbitrary Hilbert spaces $\HH$ and suited for
semi-circular systems, and the second one over an algebra
$\cC$ and suited for compound Poisson systems.

Since for the calculation of our moments we only have to deal
with elements in the algebraic Fock space (without taking a
Hilbert space completion), we will restrict ourselves to this
case in the following in order to avoid technicalities about
unbounded operators.

\subsection{Semi-circular case}

For a Hilbert space $\HH$, the algebraic full Fock space
$$\cF_{\text{alg}}(\HH):=\CC \Omega\oplus 
\HH\oplus \HH^{\otimes 2}\oplus \HH^{\otimes 3}\oplus \cdots $$
is generated by tensors $f_1\otimes\dots\otimes f_n$, where
we can think of the $f_1,\dots,f_n$ as being arranged on a
linear string.  To stress this linear nature of the usual
full Fock space, we will address it in the following as
\emph{linear Fock space} $\cF_\lin(\HH)$.  In our tracial
context, however, we should consider circular tensors, where
we think of the $f_1,\dots,f_n$ as being arranged around a
circle.  We will denote these circular tensors by $[f_1\odo
  f_n]$ and the corresponding $n$-th particle space by
$\HH_\cyclic^{\otimes n}$.  If we pair two circles, then we
have the freedom of rotating them against each other, so the
canonical inner product for this situation is given as
follows.

\begin{definition}
The \emph{cyclic Fock space} is the algebraic direct sum
\begin{equation}
\cF_\cyclic(\HH)=\bigoplus_{n=0}^\infty\HH_\cyclic^{\otimes n} 
\end{equation}
equipped with an inner product given by linear extension of
\begin{multline}
\la [f_1\odo f_n],[g_1,\odo g_m]\ra_\cyclic:=\\ \delta_{nm}\cdot
\sum_{k=0}^{n-1} \la f_1,g_{1+k}\ra\cdot \la f_2,g_{2+k}\ra\cdots
\la f_n,g_{n+k}\ra,
\end{multline}
where we count modulo $n$ in the indices of $g$.
\end{definition}

Note that one can also embed the full Fock space into the
cyclic one via
$$[f_1 \odo f_n]=\frac 1{\sqrt{n}}\sum_{k=1}^n f_k\otimes f_{k+1}\odo
f_{k-1}.$$

In order to write our formula for the fluctuations in terms
of moments of operators we still use the operators on the
linear Fock space, but we will in the end make things cyclic
by mapping the linear Fock space onto the cyclic one.

\begin{definition}

We consider the mapping $\cyc$ between linear and cyclic
Fock space,
$$\cyc:\cF_\lin(\HH)\to \cF_\cyclic(\HH),$$
which is given recursively by
\begin{equation}
\cyc \Omega=0 , \ \cyc(f) = [f],
\end{equation}
and
\begin{equation}
\cyc (f_1\otimes\cdots\otimes f_n) =[f_1 \odo f_n]
+\la f_1,\bar f_n\ra \cdot
\cyc (f_2\otimes\cdots\otimes f_{n-1})
\end{equation}
\end{definition} 

\bigskip

\begin{center}
\leavevmode
\hbox to 295pt{%
$f_1 \otimes f_2 \otimes f_3 
\otimes f_4 \otimes f_5 \sim \mbox{}
\vcenter{\hsize166pt\includegraphics{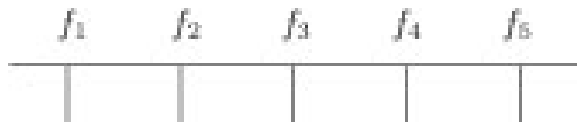}}$}

{\small\bf Figure \figno.}
\end{center}

\

\noindent
\vbox{\hsize119pt
\begin{center}
\includegraphics{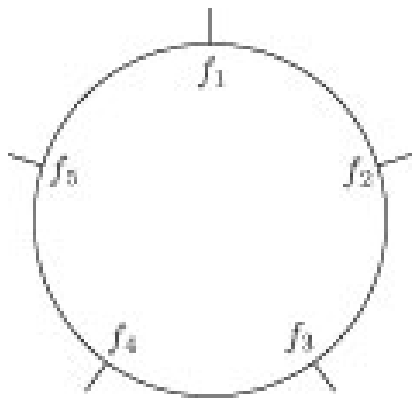}

\vskip0.5em

$
[f_1 \otimes f_2 \otimes f_3 
\otimes f_4 \otimes f_5 ]
$

\vskip0.5em
{\small\bf Figure \figno}
\end{center}}\hfill 
\vbox{\hsize119pt
\begin{center}
\includegraphics{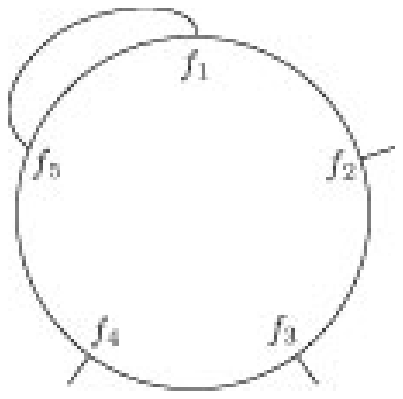}

\vskip0.5em
$
\langle f_1 , \overline{f_5} \rangle\ 
[f_2 \otimes f_3 \otimes f_4 ]
$

\vskip0.5em
{\small\bf Figure \figno}\end{center}}\hfill 
\vbox{\hsize119pt
\begin{center}

\ \includegraphics{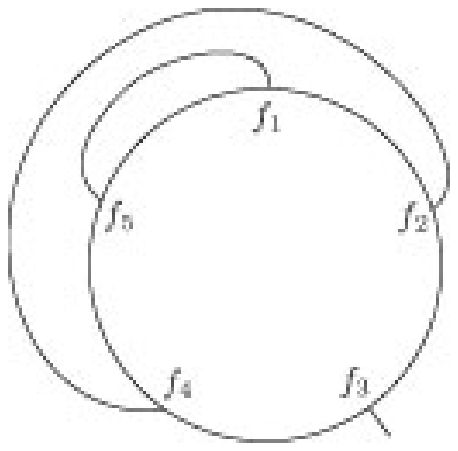}

\vskip0.5em
$
\langle f_1 , \overline{f_5} \rangle\
\langle f_2 , \overline{f_4} \rangle\ [f_3]
$
\vskip0.5em
{\small\bf Figure \figno}\end{center}}

{\leftskip=\parindent\rightskip=\parindent\noindent\small
{\bf Illustration of equation (15).} Elements of the full Fock
space, $f_1 \otimes f_2 \otimes f_3 \otimes f_4 \otimes
f_5$ are represented by linear ``half pairings'' (figure 2).
The operator {\bf c} takes a linear half pairing and wraps
it around into a circle (figure 3). Then {\bf c} pairs off
the $f$'s until either one or none remains (figures 4 and
5). The idea of a half pairing is a special case of a
half permutation explained fully in [KMS]. \par}

\bigskip

Of course, one can also write down this explicitly, here are
just two examples:
\begin{multline}
\cyc (f_1\otimes f_2\otimes f_3\otimes f_4\otimes f_5) =
[f_1\otimes f_2\otimes f_3\otimes f_4\otimes f_5]\\
+\la f_1,\bar f_5\ra \cdot[f_2\otimes f_3\otimes f_4]
+\la f_1,\bar f_5\ra\cdot \la f_2,\bar f_4\ra\cdot [f_3]
\end{multline}
and
\begin{multline}
\cyc (f_1\otimes f_2 \otimes f_3\otimes f_4\otimes f_5\otimes
f_6) = [f_1\otimes f_2\otimes f_3\otimes f_4\otimes f_5
\otimes f_6]\\+\la f_1,\bar f_6\ra \cdot[f_2\otimes f_3
\otimes f_4\otimes f_5]
+\la f_1,\bar f_6\ra\cdot \la f_2,\bar f_5\ra\cdot [f_3\otimes f_4]
\end{multline}

Let us now consider the relation between this cyclic Fock
space and fluctuations of Gaussian random matrices.  So for
the following, let $X_N(f)$ be our Gaussian random matrices
which converge, for $N\to\infty$, in distribution to a
semi-circular family, given by $\omega(f):=l(f)+l^*(f)$
realized on the full Fock space.

Our main point is now that we can express the fluctuations
of the Gaussian matrices via the operators $\omega(f)$.

\begin{theorem}\label{main1}

Let $X_N(f)$ ($f\in\HH_\RR$) be a family of Hermitian
Gaussian random matrices.  Then for all $n,m\in\NN$ and all
$f_1,\dots,f_n,g_1,\dots,\ab g_m\in\HH_\RR$
\begin{multline}\label{eq:main1}
\lim_{N\to\infty}\kk_2\{\Tr[X_N(f_1)\cdots X_N(f_n)],
\Tr[X_N(g_1)\cdots X_N(g_m)]\}\\= \la \cyc\omega(f_1) \cdots
\omega(f_n)\Omega, \cyc\omega(g_m)\cdots
\omega(g_1)\Omega\ra_{\text{\rm cyc}} \end{multline}
\end{theorem}

Note that the inversion of indices in the $g$'s is forced
upon us by the fact that our expression in random matrices
is linear in both its traces, whereas our cyclic Fock space
inner product is anti-linear in the second argument

\begin{remark}

One might wonder whether the right-hand side of our
Eq. (\ref{eq:main1}) should not also have the structure of a
variance.  This is indeed the case, but is somehow hidden in
our definition that $\cyc \Omega=0$.  If, instead of $\cyc$,
we use the mapping $\tilde \cyc$, given as follows
$$\tilde \cyc \eta:=\cyc \eta+\la \eta,\Omega\ra\Omega,$$
then the right-hand sided of $(\ref{eq:main1})$ has the
form
$$\la \cyc\eta_1,\cyc\eta_2\ra_\cyclic=
\la\tilde\cyc \eta_1,\tilde \cyc \eta_2\ra_\cyclic-
\la \tilde\cyc\eta_1,\Omega\ra_\cyclic\cdot\la \Omega,
\tilde\cyc\eta_2\ra_\cyclic.$$
\end{remark}

We will prove Theorem \ref{main1} later as a corollary of a
general limit theorem. For the moment, we will be content
with checking the consistency of our statement with respect
to traciality.  Since the left hand side is tracial in the
arguments of the traces, the right hand side should be
tracial, too.  Recall that $\cA(\HH_\RR)$ is the unital
$*$-algebra generated by all $\omega(f)=l(f)+l^*(\bar f)$
with $f\in\HH$.

\begin{lemma}
The mapping
\begin{align}
\cA(\HH_\RR)&\to\cF_\cyclic(\HH)\\
a&\mapsto \cyc a\Omega
\end{align}
is tracial, i.e., for all $a,b\in \cA(\HH_\RR)$ we have
$$\cyc ab\Omega=\cyc ba\Omega.$$
\end{lemma}

\begin{proof}

Since $\Omega$ is cyclic and separating for the algebra
$\cA(\HH_\RR)$, it suffices to show that
$$\cyc \omega(f)W(f_1\odo f_n)\Omega=\cyc W(f_1\odo f_n)\omega(f)
\Omega$$

\noindent
for all $f,f_1,\dots,f_n\in\HH$.  On the left side, we have
$$
\cyc \omega(f)W(f_1\odo f_n)\Omega=
\cyc (f\otimes f_1\odo f_n) +
\la f_1,\bar f\ra \cdot  \cyc (f_2\odo f_n).
$$
For the right side, it follows from the relation\footnote{As
we have been unable to locate a proof of this in the
literature one is given in \cite[\S 10]{KMS}. The main idea is to write
$\omega(f_1) \omega(f_2) \cdots \omega(f_n)$ as a linear combination of
Wick polynomials. In fact we may write $\omega(f_1) \omega(f_2)
\cdots
\omega(f_n) = \sum_\pi W_\pi(f_1 \otimes f_2 \otimes \cdots \otimes f_n)$
(*) where the sum is over all non-crossing ``half-pairings'' of $[n]$,
that is non-crossing pairings of $[n]$ with only singletons and pairs
such that the singletons are not enclosed by any pair, and $W_\pi(f_1
\otimes \cdots \otimes f_n) = \langle f_{i_1}, \overline{f_{j_1}}\rangle
\cdots
\langle f_{i_k},
\overline{f_{j_k}}\rangle W(f_{l_1} \otimes \cdots \otimes f_{l_m})$ where
the pairs of
$\pi$ are $(i_1, j_1), \dots, (i_k, j_k)$ and the singletons are $(l_1),
\dots, (l_m)$. Note that in (*) there is only one term involving a Wick
polynomial on $n$ $f$'s. So that if we know that $W(f_1 \otimes \cdots
\otimes f_m)^\ast = W(\overline{f_m} \otimes \cdots \otimes
\overline{f_1})$ for $m < n$ they we apply induction to the adjoint of
(*). }\label{firstfootnote}
$W(f_1 \otimes \cdots \otimes f_n)^\ast = W(\bar f_n \otimes \cdots
\otimes
\bar f_1)$ that
$$
W(f_1\odo f_n)\omega(f)=W(f_1\odo f_n\otimes f)+\la f_n,\bar f\ra
W(f_1\odo f_{n-1}),
$$ 
which yields
$$
\cyc W(f_1 \odo f_n) \omega(f) \Omega= \cyc (f_1 \odo f_n
\otimes f) + \la f_n, \bar f \ra \cdot \cyc ( f_1 \odo
f_{n-1}).
$$ 
From the definition of $\cyc$ we see that both sides are the
same.
\end{proof}

\subsection{Compound Poisson case}

Let us now consider the case where we have a $*$-algebra
$\cC$ with trace $\psi$.  We denote by
$$\cF_\lin(\cC)=\bigoplus_{n=0}^\infty \cC^{\otimes n}$$ 
the algebraic linear Fock space and by
$$\cF_{\cyclic}(\cC)=\bigoplus_{n=0}^\infty \cC_\cyclic^{\otimes n}$$
the algebraic cyclic Fock space.

Since in this case we also have actions of our operators
which multiply inside the argument, we have to take this
into account when we glue the beginning and end of the
tensors together. Thus we have to change the definition of
the map $\cyc$ as follows.

\begin{definition}
We consider the linear mapping
$$\cyc:\cF_\lin(\cC)\to \cF_\cyclic(\cC),$$
given recursively by
$$\cyc\Omega:=0 \mbox{\ and\ } \cyc(d) = [d]$$
and
\begin{multline}\label{poissoncdef}
\cyc (d_1\otimes\cdots\otimes d_n) := [d_1 \odo d_n]+
[d_nd_1\otimes d_2 \odo d_{n-1}]\\+
\psi(d_1d_n)\cdot \cyc (d_2\otimes \dots\otimes d_{n-1} )
\end{multline}
\end{definition}

For example
\begin{equation}\label{smallexample}
\cyc(d_1 \otimes d_2 \otimes d_3) = [d_1 \otimes d_2 \otimes d_3] + [d_3d_2 \otimes
d_2] + \psi(d_3d_1) [d_2]
\end{equation}

\bigskip

\begin{center}
$d_1 \otimes d_2 \otimes d_3 \sim
\vcenter{\hsize110pt\includegraphics{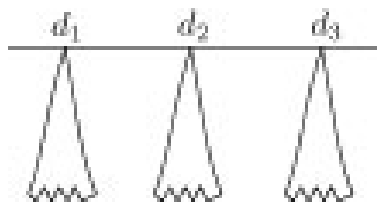}}$\\[10pt]
{\small \bf Figure \figno.}
\end{center}

\leavevmode\kern-20pt
\vbox{\hsize125pt\begin{center}\includegraphics{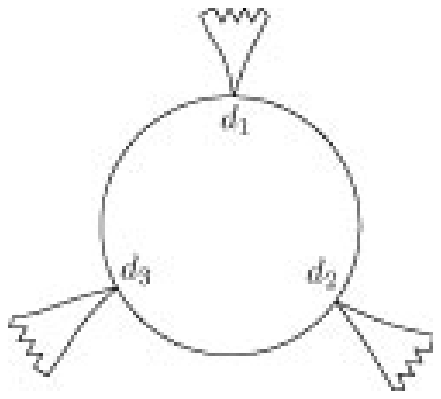}\\
$[d_1 \otimes d_2 \otimes d_3]$\\[10pt]
{\small\bf Figure \figno.}\end{center}}
\vbox{\hsize125pt\begin{center}\includegraphics{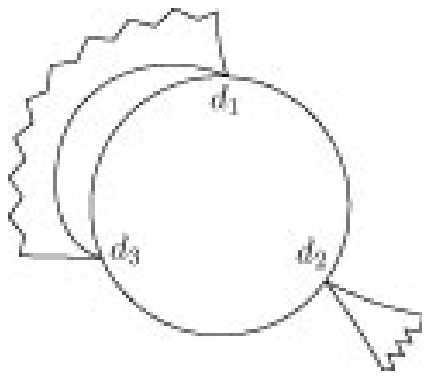}\\
$[d_3d_1 \otimes d_2]$\\[10pt]
{\small\bf Figure \figno.}\end{center}}
\vbox{\hsize125pt\begin{center}\includegraphics{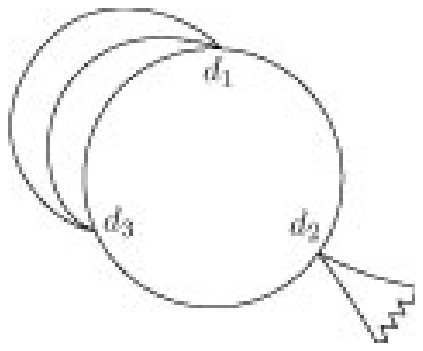}\\
$\psi(d_3d_1) [d_2]$\\[10pt]
{\small\bf Figure \figno.}\end{center}}

{\leftskip=\parindent\rightskip=\parindent\noindent\small
{\bf Illustration of Equation (\ref{smallexample})} The
vector $d_1 \otimes d_2 \otimes d_3$ is represented by a
linear half permutation (figure 6) with one open block for
each factor in the tensor product ({\it c.f.} \cite[\S
10]{KMS}). $[d_1 \otimes d_2 \otimes d_3]$ is represented by
a circular half permutation with one open block for each
factor in the tensor product. $\psi(d_3d_1) [d_2]$ is
represented by a circular half permutation with one closed
block and one open block (figure 9). The operator $\cyc$
first turns the linear half permutation into a circular half
permutation (figure 7). Then $\cyc$ fuses a pair of open
blocks (figure 8) and then closes the just formed open block
(figure 9). This process continues until either one or zero
open blocks remain.  \par}

\bigskip\bigskip

Then we claim that one can express the fluctuations of our
compound Wishart matrices also by calculations in terms of
the corresponding operators $p(d)$.

\begin{theorem}\label{main2}

Suppose $\{D^{(N)}_1, \dots , D^{(N)}_p\}$ converges in
distribution to $\{d_1, \dots , d_p\}$ in $(\cC, \psi)$ and
$(X_N)_{N\in\NN}$ is a sequence of $N\times N$ Hermitian
Gaussian random matrices. We put $P_N(D_i):=X_ND^{(N)}_iX_N$ and
let $p(d_i)$ be our operators on the full Fock space, then we
have for all $m,n\in\NN$ that
\begin{multline}
\lim_{N\to\infty}
\kk_2\{\Tr[P_N(D_{i_1})\cdots P_N(D_{i_n})],
\overline{\Tr[P_N(D_{i_{n+1}})\cdots  P_N(D_{i_{n+m}})] }\} \\
= \la \cyc p(d_{i_1}) \cdots p(d_{i_n}) \Omega, \cyc
p( d_{i_{n+1}} )\cdots p(d_{i_{n+m}}) \Omega\ra_\cyclic.
\end{multline}
\end{theorem}

Again, we check only the traciality of the right hand side
and postpone the proof of the statement until we have proved
our general limit theorem.  Recall that we denote by
$\cA(\cC)$ the unital $*$-algebra generated by all $p(d)$
for $d\in\cC$.

\begin{lemma}\label{tracial-poisson}

The mapping
\begin{align*}
\cA(\cC)&\to \cF_\cyclic(\cC)\\
a&\mapsto \cyc a\Omega
\end{align*}
is tracial.
\end{lemma}

\begin{proof}
Since $\Omega$ is cyclic and separating for $\cA(\cC)$
(see \cite{Ans}), it suffices to check for $d,d_1,\dots,d_n\in\cC$
that $$\cyc p(d)W(d_1\odo d_n)\Omega =\cyc W(d_1\odo d_n)p(d)\Omega.$$
For $n=0$, i.e., $W(\Omega)=1$, this is surely true. In general we
have for the left hand side
\begin{align*}
\cyc p(d)W(d_1&\odo d_n)\Omega=
\cyc p(d)d_1\odo d_n\\
&=\cyc\bigl( d\otimes d_1\odo d_n+ \psi(d d_1) d_2\odo d_n \\
&\qquad+
dd_1\otimes d_2\odo d_n+\psi(d) d_1\odo d_n\bigr)
\end{align*}
By using the identity\footnote{%
The proof is very similar to that
sketched in the footnote on page~\pageref{firstfootnote}.
A detailed proof is provided in \cite[\S 10]{KMS}.}
$W(d_1
\otimes
\cdots
\otimes d_n)^\ast = W(d^\ast_n \otimes \cdots \otimes d_1^\ast)$ we have
\begin{align*}
W(d_1\odo d_n)p(d)&=W(d_1\odo d_n\otimes d)\\&\qquad+ \psi(d_nd)W(d_1\odo
d_{n-1})\\
&\qquad+W(d_1\odo d_{n-1}\otimes d_nd)\\&\qquad+
\psi(d) W(d_1\odo d_n),
\end{align*}
Thus the right hand side becomes
\begin{align*}
\cyc W(d_1\odo d_n)p(d)\Omega
&=\cyc\bigl( 
d_1\odo d_n\otimes d \\&\qquad+\psi(d_nd)d_1\odo
d_{n-1}\\
&\qquad+d_1\odo d_{n-1}\otimes d_nd\\&\qquad+
\psi(d) d_1\odo d_n
\bigr)
\end{align*}
So it remains to show that
\begin{multline*}
\cyc\bigl( d\otimes d_1\odo d_n+ \psi(d d_1) d_2\odo d_n +
dd_1\otimes d_2\odo d_n\bigr)\\
=\cyc\bigl( 
d_1\odo d_n\otimes d+ \psi(d_nd)d_1\odo
d_{n-1}+d_1\odo d_{n-1}\otimes d_nd
\bigr)
\end{multline*}
This can be checked directly by applying the definition of the
mapping $\cyc$.
\end{proof}

\section{Second order freeness and abstract limit theorems}

We shall derive our main theorems, \ref{main1} and
\ref{main2}, from a general limit theorem, very much in
the same spirit as one can get the distribution of the
semi-circle and the compound free Poisson distributions from
free limit theorems, see \cite{Sp1}. The crucial idea is the
notion of second order freeness which we introduce below.

\begin{definition}
A \emph{second order non-commutative probability space} \ab
$(\cA,\ab    \ff,\nappa)$ consists of a unital algebra $\cA$,
a tracial linear functional 
$$\ff:\cA\to\CC \qquad\text{with} \qquad \ff(1)=1$$ 
and a bilinear functional 
$$\nappa:\cA\times\cA\to\CC,$$ 
which is tracial in both arguments and which satisfies
$$\nappa(a,1)=0=\nappa(1,b)\qquad\text{for all $a,b\in\cA$.}$$
\end{definition}

\begin{notation}

Let unital subalgebras $\cA_1,\dots,\cA_r\subset\cA$ be given.

\noindent
1) We say that a tuple $(a_1,\dots,a_n)$ ($n\geq 1$) of
elements from $\cA$ is \emph{cyclically alternating} if, for
each $k$, we have an $i(k)\in\{1,\dots,r\}$ such that
$a_k\in \cA_{i(k)}$ and, if $n\geq 2$, we have $i(k)\not=
i(k+1)$ for all $k=1,\dots,n$. We count indices in a cyclic
way modulo $n$, i.e., for $k=n$ the above means
$i(n)\not=i(1)$.  Note that for $n=1$ we mean that $a_1$ is
in some $\cA_i$.  

\noindent
2) We say that a tuple $(a_1,\dots,a_n)$ of elements from
$\cA$ is \emph{centered} if we have
$$\ff(a_k)=0\qquad\text{for all $k=1,\dots,n$.}$$
\end{notation}

\begin{definition}\label{free}

Let $(\cA,\ff,\nappa)$ be a second order non-commutative
probability space.  We say that unital subalgebras
$\cA_1,\dots,\cA_r\subset \cA$ are \emph{free with respect
to $(\ff,\nappa)$} or \emph{free of second order}, if they
are free with respect to $\ff$ and whenever we have centered
and cyclically alternating tuples $(a_1,\dots,a_n)$ and
$(b_1,\dots, b_m)$ from $\cA$ then we have:
\begin{enumerate}
\item $\nappa(a_1\cdots a_n, b_1\cdots b_m)=0$ for $n\not=m$;
\item $\nappa(a,b)=0$ for $a\in\cA_i$, $b\in\cA_j$, and $i\not= j$;
\item if $n=m\geq2$, then
$$\nappa(a_1\cdots a_n, b_1\cdots b_n)=
\sum_{k=0}^{n-1} \ff(a_1b_{n+k})\cdot\ff(a_2b_{(n-1)+k})\cdots
\ff(a_nb_{1+k}).$$
\end{enumerate}

Note that in the sum the indices of the $a$'s increase,
whereas those of the $b$'s decrease; one should think of two
concentric circles with the $a$'s on one of them and the
$b$'s on the other. However, whereas on one circle we have a
clockwise orientation of the points, on the other circle the
orientation is counter-clockwise. Thus, in order to match up
these points modulo a rotation of the circles, we have to
pair the indices as in the sum above.
\end{definition}

Condition ({\it iii}) is the annular version of the disc
picture of first order freeness: suppose $a_k \in \cA_{i_k}$ with
$\phi(a_k) = 0$ for $1 \leq k \leq n$ and we arrange the
elements $a_1, \dots, a_n$ around the boundary of a disc
(figure 10). The only non-crossing partition of $[n]$ that
only connects elements from the same algebra consists of all
singletons and since $\phi$ of a singleton is 0, we have
that $\phi(a_1 \cdots a_n) = 0$.

In the annular case we put the centered and cyclically
alternating elements $(a_1, \dots, a_n)$ and $(b_1, \dots,
b_n)$ around the boundary of an annulus. We only connect
elements from the same algebra and as the elements are
centered we have no singletons; so we must connect in pairs
elements from opposite circles is all possible ways (figure
11). This is the meaning of condition ({\it iii}).

\begin{center}\leavevmode
\vbox{\hsize125pt\begin{center}\includegraphics{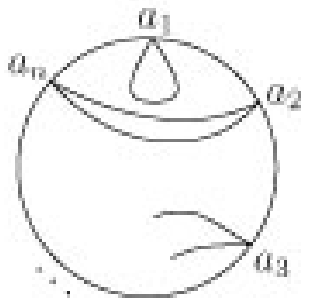}\\
{\small\bf Figure \figno.}\end{center}}\qquad
\vbox{\hsize125pt\begin{center}\includegraphics{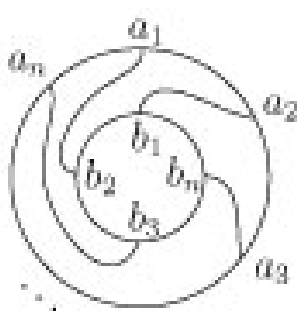}\\
{\small\bf Figure \figno.}\end{center}}
\end{center}

Note that, as in the case of freeness, the trick of writing
elements $a$ as 
$$a=a^{\oh}+\ff(a)\cdot 1,\qquad\text{where $\ff(a^{\oh})=0$,}$$
allows us to calculate $\nappa$ in terms of $\ff$ and
$\nappa$ restricted to the subalgebras. However, whereas the
formulas for $\ff$ of mixed moments contain only $\ff$
applied to the subalgebras, $\nappa$ of mixed moments has in
general to be expressed in both $\ff$ and $\nappa$
restricted to the subalgebras.

For example, assume we have two subalgebras $\cA_1$ and
$\cA_2$, and elements $a_1,a_2\in\cA_1$ and
$b_1,b_2\in\cA_2$. Then we have
\begin{equation}\label{e24}
\nappa(a_1,b_1)=0,
\end{equation}
\begin{equation}\label{e25}
\nappa(a_1 b_2, a_2)=\ff(b_2)\cdot \nappa(a_1,a_2)
\end{equation}
or
\begin{multline}\label{beispiel}
\nappa(a_1b_2,a_2b_2) =
\ff(a_1a_2)\ff(b_1b_2)-\ff(a_1a_2)\ff(b_1) \ff(b_2) \\
\mbox{} -\ff(a_1)\ff(a_2)\ff(b_1b_2) +\ff(a_1)\ff(a_2)\ff(b_1)
\ff(b_2) \\
\mbox{} +\nappa(a_1,a_2)\ff(b_1)\ff(b_2) +
\ff(a_1)\ff(a_2)\nappa(b_1,b_2).
\end{multline}

One should note that each variable appearing in the
arguments of $\nappa$ on the left-hand side of these
examples has to appear exactly once in each product on the
right-hand side. Let us formalize this in the following
definition.
 
\begin{notation}

Let $(\cA,\ff,\nappa)$ be a second order non-commutative
probability space with subalgebras $\cA_1,\dots,\cA_r \subset \cA$,
and consider elements $a_1,\dots,a_n\in \cup_{i=1}^r\cA_i$. A
\emph{balanced expression} (with respect to the subalgebras $\cA_1,
\dots , \cA_r$)  in $a_1,\dots,a_n$ is a product of factors
$\ff(a_{i_1} \cdots a_{i_t})$ and $\nappa(a_{i_1} \cdots
a_{i_s},\ab a_{j_1} \ab\cdots a_{j_t})$ where each $a_i$ has to
appear exactly once among all arguments and the argument of each
$\ff$ or the arguments of each $\nappa$ contains only $a_i$'s from
a single $\cA_j$.

\end{notation}

For example, balanced expressions in $a_1,a_2,a_3,a_4$ are
$$
\ff(a_1a_3)\ff(a_2a_4) \mbox{ if } a_1, a_3 \in \cA_1 \mbox{
and } a_2, a_4 \in \cA_2$$ $$\text{or}\quad
\ff(a_1)\ff(a_4)\nappa(a_2,a_3) \mbox{ if } a_1, a_4 \in
\cA_1 \mbox{ and } a_2, a_3 \in \cA_2
$$ 
Every summand on the right-hand side of Eq. (\ref{beispiel})
is a balanced expression in $a_1,a_2,b_1,b_2$ if $a_1, a_2
\in \cA_1$ and $b_1, b_2 \in \cA_2$.

\begin{lemma}\label{fact:good}
Let $\cA_1,\dots,\cA_r$ in $(\cA,\ff,\nappa)$ be free with
respect to $(\ff,\nappa)$. Suppose we have cyclically
alternating $(a_1,\dots,a_n)$ and $(b_1,\dots,b_m)$ and
denote by $s$ the number of different subalgebras appearing
in $\{a_1, \dots a_n,\ab b_1, \dots b_m\}$.  Then
$\nappa(a_1\cdots a_n, b_1\cdots b_m)$ is either 0 or can be
written as a sum of balanced expressions in
$a_1,\dots,a_n,b_1,\dots,b_m$, such that each of these
balanced expressions has at least $s$ factors and contains
at most one $\nappa$-factor.

Thus any expression of the form $\nappa(a_1\cdots a_n,
b_1\cdots b_m)$ for cyclically alternating $(a_1,\dots,a_n)$
and $(b_1,\dots,b_m)$ is determined by the value of $\ff$
restricted to $\cA_1\cup\dots \cup \cA_r$ and by the value
of $\nappa$ restricted to $(\cA_1\times \cA_1)\cup \dots\cup
(\cA_r\times \cA_r)$.
\end{lemma}

\begin{proof}
We will prove this by induction on $n+m$. The case $n+m=2$,
i.e., $n=m=1$, is clear.

So consider $n+m\geq 3$.  Put
$$
a_k^{\oh}:=a_k-\ff(a_k)\cdot 1,\qquad
b_l^{\oh}:=b_l-\ff(b_l)\cdot 1
$$
for $k=1,\dots,n$ and $l=1,\dots,m$.  Then we have
\begin{align}\label{e27}
&\nappa(a_1\cdots a_n, b_1\cdots b_m) \nonumber\\
& = \nappa\bigl((a_1^{\oh}+\ff(a_1))\cdots (a_n^{\oh}+\ff(a_n)), 
(b_1^{\oh}+\ff(b_1))\cdots (b_m^{\oh}+\ff(b_m))\bigr)\nonumber\\
&=\sum_{p,q}
\ff(a_{p(1)})\cdots \ff(a_{p(k)})\cdot\ff(b_{q(1)}) \nonumber\\
& \hbox{}\kern70pt \hfill\times \cdots \times
\ff(b_{q(l)})\cdot \nappa\bigl(a_{\bar p(1)}^{\oh}\cdots
a_{\bar p(n-k)}^{\oh},b_{\bar q(1)}^{\oh}\cdots b_{\bar q(m-l)}^{\oh}\bigr),
\end{align}
where the sum runs over all partitions 
$$
((p(1),\dots,p(k)), (\bar p(1),\dots,\bar
p(n-k))\qquad\text{of the set $[1,n]$}
$$
and
$$
((q(1),\dots,q(l)), (\bar q(1),\dots,\bar
q(m-l))\qquad\text{of the set $[1,m]$}
$$ 
into two ordered subsets (with $k=0,\dots,n$ and
$l=0,\dots,m$).  The term corresponding to $k = l =0$ is, by
Definition \ref{free}, either 0 (when $m \not= n$) or is a
balanced expression in the centered elements with at least
one factor for each occurring subalgebra. Now note that a
balanced expression in centered elements can be rewritten as
a sum of balanced expressions in the original elements and
that the number of factors can only increase by doing so.

For the other terms with $k+l\geq 1$, $(a^{\oh}_{\bar p(1)},
\dots , a^{\oh}_{\bar p(n-k)})$ and $(b^{\oh}_{\bar q(1)},
\dots ,\ab b^{\oh}_{\bar q(m-l)})$ may no longer be
cyclically alternating. So we group together adjacent
elements from the same algebra to produce a cyclically
alternating tuple with at least $\max \{1,s - (k+l)\}$ subalgebras
appearing, and so we can apply our induction hypothesis.
Indeed, the term
\begin{equation}\label{kasteq}
\nappa\bigl(a_{\bar p(1)}^{\oh}\cdots
a_{\bar p(n-k)}^{\oh},b_{\bar q(1)}^{\oh}\cdots b_{\bar q(m-l)}^{\oh}\bigr)
\end{equation}
contains elements from at least $s-(k+l)$ different
subalgebras; by our induction hypothesis, we may write it as
the sum of balanced expressions in the $a^{\oh}$'s and
$b^{\oh}$'s, each product containing at least $s-(k+l)$
factors.  Again we write a balanced expression in centered
elements as a sum of balanced expressions in the original
elements.  This means we can write the term (\ref{kasteq})
as a sum of balanced expressions in
$$
a_{\bar p(1)}, \cdots, a_{\bar p(n-k)},b_{\bar q(1)},\cdots,
b_{\bar q(m-l)}
$$ 
with at least $s-(k+l)$ factors for each product. Together
with the $k+l$ factors
$$
\ff(a_{p(1)})\cdots \ff(a_{p(k)})\cdot\ff(b_{q(1)})\cdots
\ff(b_{q(l)})
$$
this gives the assertion. Note that in all our steps
balancedness is preserved and that at most one $\nappa$-term
can occur in all the reductions.
\end{proof}

A very special case of such a factorization is given in the
next lemma.

\begin{lemma}\label{properties}

Let $(\cA,\ff,\nappa)$ be a second order non-commutative
probability space and let $\cA_1,\dots,\cA_r\subset\cA$ be
free with respect to $(\ff,\nappa)$.  Consider cyclically
alternating $(a_1\dots,a_n)$ and $(b_1,\dots, b_m)$ from
$\cA$.

\noindent
1) Assume that the subalgebra of $a_1$ appears only once. Then
we have
$$
\nappa(a_1\cdots a_n,b_1\cdots b_m)=\ff(a_1)\nappa(a_2\cdots
a_n,b_1\cdots b_m).
$$
2) Assume that the subalgebra of $b_1$ appears only once. Then we have
$$\nappa(a_1\cdots a_n,b_1\cdots b_m)=\ff(b_1)\nappa(a_1\cdots
a_n,b_2\cdots b_m).$$
\end{lemma}

\begin{proof}
We only prove the first part.  Put
$$
a_k^{\oh}:=a_k-\ff(a_k)1,\qquad
b_l^{\oh}:=b_l-\ff(b_l)1.
$$
We have
\begin{align*}
\nappa(a_1\cdots &a_n,b_1\cdots b_m)=\nappa \bigl((a_1^{\oh}+
\ff(a_1)1)a_2\cdots
a_n,b_1\cdots b_m\bigr)\\
&=\nappa(a_1^{\oh} a_2\cdots a_n,b_1\cdots b_m)+
\ff(a_1)\nappa(a_2\cdots a_n,b_1\cdots b_m).
\end{align*}
We shall show that the first term is 0. 

Indeed, we shall show that if $(a_1, a_2, \dots , a_n)$ and$(b_1,
b_2, \dots , b_m)$ are cyclically alternating and the algebra of
$a_1$ appears only once then $\nappa(a^{\oh}_1 a_2\ab \cdots \ab
a_n, b_1 \cdots b_m) = 0$. We shall do this by induction on $m+n$.
By equations (\ref{e24}) and (\ref{e25}) we have
$\nappa(a^{\oh}_1 a_2 \cdots \ab a_n, b_1 \cdots b_m) = 0$ when $m+n
= 2$ or 3. Suppose we have proved the result for $m + n <
j$; we shall prove it for $m + n = j$.

We shall use the expansion in equation (\ref{e27}) and show that
for $ 1 \leq k \leq n-1$ and $1 \leq l \leq m$, $\nappa(a_1^{\oh}
a_{p(1)}^{\oh} \cdots a_{p(k)}^{\oh}, b^{\oh}_{q(1)} \cdots
b^{\oh}_{q(l)}) = 0$ for all subsets $\{ p(1), \dots , p(k)\}
\subset \{1, 2, 3, \dots , n-1\}$ and $\{q(1), \dots , q(l)\}
\subset \{1, 2, 3, \dots , m\}$. 

When $k= n-1$ and $l=m$ we have that $(a_1^{\oh}, \dots ,
a_n^{\oh})$ and $(b^{\oh}_1, \dots , b^{\oh}_m)$ are
centered and cyclically alternating. If $m \not = n$ we have
$\nappa(a^{\oh}_1 a_2^{\oh} \cdots \ab a_n^{\oh}, b_1^{\oh} \cdots
b_m^{\oh}) = 0$ by {\it (i)} of Definition \ref{free}. If $m = n$
then by {\it(iii)}, we have $\nappa(a^{\oh}_1 a_2^{\oh} \cdots \ab
a_n^{\oh}, b_1^{\oh} \cdots b_m^{\oh}) = 0$ because $\ff(a_1^{\oh}
b_i^{\oh}) = 0$ for all $i$.

Suppose next that $k + l \leq m + n -2$. We can no longer
expect $(a_1^{\oh}, a_{p(1)}^{\oh}, \dots , a_{p(k)}^{\oh})$
and $(b_{q(1)}^{\oh}, \cdots , b_{q(l)}^{\oh})$ to be
cyclically alternating; so we group adjacent terms from the
same algebra and write $a_1^{\oh} a_{p(1)}^{\oh} \cdots\ab
a_{p(k)}^{\oh}\ab = a_1^{\oh} c_1 \cdots c_s$ and
$b_{q(1)}^{\oh}\cdots b_{q(l)}^{\oh} = d_1 \cdots d_t$ with
$(c_1, \dots, c_s)$ and $(d_1, \dots , d_t)$ cyclically
alternating and neither involving the algebra of $a_1$. Now
$s + t + 1 \leq k +l + 1 \leq m + n -1 \leq j-1$, so by our
induction hypothesis $\nappa(a^{\oh}_1 a_{p(1)}^{\oh} \cdots \ab
a_{p(k)}^{\oh}, b_{q(1)}^{\oh} \cdots b_{q(l)}^{\oh}) =
\nappa(a_1^{\oh} c_1 \cdots c_s, d_1 \cdots  d_t) = 0$
\end{proof}

If by the lemma above we successively remove all subalgebras
which occur only once and multiply together cyclic
neighbours from the same subalgebra, then we arrive finally
at $\nappa(a,b)$ for $a,b$ from one of the subalgebras (both
from the same, in order to get a non-vanishing contribution)
or at $\nappa(a_1\cdots a_n, b_1\cdots b_m)$ where both
arguments are cyclically alternating and in addition each
involved subalgebra appears at least twice. In the latter
case we have either a very special matching of the involved
subalgebras or we can strengthen Lemma \ref{fact:good} to
obtain at least one more $\ff$-factor.

\begin{lemma}\label{fact:bad}

Let $(\cA,\ff,\nappa)$ be a second order non-commutative
probability space and let $\cA_1,\dots,\cA_r\subset\cA$ be
free with respect to $(\ff,\nappa)$. Suppose
$(a_1,\dots,a_n)$ and $(b_1,\dots,b_m)$ are cyclically
alternating and denote by $s$ the number of different
subalgebras appearing in $\{ a_1, \dots , a_n, b_1,\ab \dots
b_m \}$. Suppose also that each involved subalgebra appears
at least twice.

Then $\nappa(a_1 \cdots a_n, b_1 \cdots b_m)$ can be written as a
sum of balanced expressions with at least $s+1$ factors unless
the following conditions are satisfied: 
$$
(*)\ \left\{\vcenter{\hsize320pt\noindent --- \ $m=n=s$;\\
--- \ for each $k$ there is $k'$ such that $a_k$ and
$b_{k'}$ are from the same subalgebra; and \\ --- \ there
is $q$ such that for all $k$, $k'=-k+q \, \mod
n$. }\right.
$$ In this case $\nappa(a_1 \cdots a_n, b_1
\cdots b_n) = \ff(a_1 b_{1'})  \cdots  \ff(a_n b_{n'}) + S$, where
$S$ is a sum of balanced expressions with at least $s+1$ factors.
\end{lemma}

\begin{proof}
Let us look again at the expansion
\begin{align}
&\nappa(a_1\cdots a_n, b_1\cdots b_m)\nonumber\\
&=\sum_{p,q}
\ff(a_{p(1)})\cdots \ff(a_{p(k)})\cdot\ff(b_{q(1)})\cdots
\ff(b_{q(l)}) \nonumber\\
&\qquad\mbox{} \hfill\cdot \nappa\bigl(a_{\bar p(1)}^{\oh}\cdots
a_{\bar p(n-k)}^{\oh},b_{\bar q(1)}^{\oh}\cdots b_{\bar q(m-l)}^{\oh}\bigr).
\end{align}

First, consider a term with $k+l \geq 1$. Then there are two
possibilities. If all $\{a_{p(1)},\dots,a_{p(k)},
b_{q(1)},\dots, b_{q(l)} \}$ belong to different
subalgebras, then there must be exactly $s$ subalgebras in
$\{a_{\bar p(1)}, \dots , a_{\bar p(n-k)}, b_{\bar q(1)}
,\ab \dots ,\ab b_{\bar q(m-l)} \}$ because each involved
subalgebra appears at least twice. If we group together any
adjacent terms that may come from the same subalgebra we
obtain cyclically alternating arguments and so can apply
Lemma \ref{fact:good}.  According to Lemma \ref{fact:good}
we can write $\nappa(a_{\bar p(1)}^{\oh} \cdots a_{\bar
p(n-k)}^{\oh},\ab b_{\bar q(1)}^{\oh}\ab \cdots b_{\bar
q(m-l)}^{\oh})$ as a sum of balanced expressions with at
least $s$ factors. Combining these with the $k +l $ factors
$\ff(a_{p(1)}) \cdots \ff(a_{p(k)}) \ff(b_{q(1)}) \cdots \ab
\ff(b_{q(l)})$ we have that every term with $k+l \geq 1$ can
be written as a sum of balanced expressions with at least
$s+1$ factors.

Second, consider the term $\nappa(a_1^{\oh} \cdots
a_n^{\oh}, b_1^{\oh} \cdots b_m^{\oh})$ corresponding to $k
= l = 0$. If $m \not = n$ we have this is zero by Definition
\ref{free}. So suppose $m = n$. Again by Definition
\ref{free} $\nappa(a_1^{\oh} \cdots a_n^{\oh}, b_1^{\oh}
\cdots b_n^{\oh}) = \sum_{k=0}^{n-1} \ff(a^{\oh}_1
b^{\oh}_{n+k}) \cdots\ab \ff(a^{\oh}_n b^{\oh}_{k+1})$.

If $s < n$ then each term $\ff(a^{\oh}_1 b^{\oh}_{n+k})
\cdots\ab \ff(a^{\oh}_n b^{\oh}_{k+1})$ has $n \geq s+1$
factors or is zero if for some factor $\ff(a^{\oh}_r
b^{\oh}_{n+k-r})$, $a_r$ and $b_{n+k-r}$ come from different
algebras. Thus we get either 0 or a balanced expression with
at least $s+1$ factors.

Finally assume that $s=m=n$. Each subalgebra must appear
exactly twice, so for each $k$ there is $k'$ such that $a_k$
and $b_{k'}$ are from the same subalgebra, or else for all
$q$, $\ff(a^{\oh}_1 b^{\oh}_{n+q - 1 }) \cdots\ab \ff(a^{\oh}_n
b^{\oh}_{q}) = 0$. Again we will have $\ff(a^{\oh}_1
b^{\oh}_{n+q-1}) \cdots\ab \ff(a^{\oh}_n b^{\oh}_{q}) = 0$
unless $k' = -k + q \mod n$, for some $q$. For this $q$ we have
$\nappa(a_1^{\oh} \cdots a_n^{\oh}, b_1^{\oh} \cdots
b_n^{\oh}) = \ff(a^{\oh}_1 b^{\oh}_{1'}) \cdots\ab
\ff(a^{\oh}_n b^{\oh}_{n'})$.

By substituting $\ff(a^{\oh}_k b^{\oh}_{k'}) = \ff(a_k
b_{k'}) - \ff(a_k) \ff(b_{k'})$ into $
\ff(a^{\oh}_1 b^{\oh}_{1'}) \cdots\ab \ff(a^{\oh}_n
b^{\oh}_{n'})$ we may write $\nappa(a_1^{\oh} \cdots
a_n^{\oh}, b_1^{\oh} \cdots b_n^{\oh})$ as \ab $\ff(a_1
b_{1'}) \cdots\ab \ff(a_n b_{n'})$ plus a sum of balanced
expressions with at least $n+1 = s+1$ factors.
\end{proof}

We are now almost ready for the main limit theorem of
second order freeness. It will turn out that moments of the
limit can be calculated in terms of annular non-crossing
objects. However, in this setting we will not arrive
directly at permutations (as in the fluctuation formulas for
random matrices), but -- as is much more natural in the
context of limit theorems -- at partitions.  In the random
matrix setting of section 3 we got contributions of the form
$\psi_\pi$ for non-crossing permutations $\pi$. So we have
to define the analogous object $\csi_\sigma$ for
non-crossing partitions $\sigma$.  However, for non-crossing
partitions, the contribution to $\csi_\sigma$ of a block
which is the only through-block will require special
treatment. We will need two different types of functions,
$\csi_1$ in the case of multiple through-blocks and $\csi_2$
in the case of a single through-block.

\begin{notation}\label{not:oje}
Let $\cI$ be an index set and let two functions 
\[
\csi_1:\bigcup_{n\in\NN} \cI^n \to\CC, \qquad
(t_1,\dots,t_n)\mapsto \csi_1(t_1,\dots,t_n)
\]
and 
$$\displaylines{
\csi_2:\bigcup_{n,m\in\NN} \cI^n\times \cI^m \to\CC \cr
(t_1,\dots,t_n) \times (t_{n+1}, \dots , t_{n + m}) \mapsto
\csi_2(t_1,\dots,t_n; t_{n+1} , \dots , t_{n + m}) \cr}$$
be given.
Assume that $\csi_1$ is tracial in its arguments, i.e.,
for all $n\in\NN$ and all $t_1,\dots,t_n\in \cI$ we have
$$\csi_1(t_1,t_2,\dots,t_n)=\csi_1(t_2,\dots,t_n,t_1),$$
and that $\csi_2$ is tracial in each of its groups of arguments, i.e.,
for all $n,m\in\NN$ and all $t_1,\dots,t_{n+m}$ we have
$$\csi_2(t_1,\dots,t_n;t_{n+1},\dots,t_{n+m})=
\csi_2(t_2,\dots,t_n,t_1;t_{n+1},\dots,t_{n+m})$$
and
$$\csi_2(t_1,\dots,t_n;t_{n+1},\dots,t_{n+m})=
\csi_2(t_1,\dots,t_n;t_{n+2},\dots,t_{n+m},t_{n+1}).$$

Fix $n,m\in\NN$ and consider an annular non-crossing
partition $\sigma\in \NC(n,m)$.  Then, for given
$t_1,\dots,t_n,t_{n+1},\dots t_{n+m}\in\cI$ we define
$\csi_\sigma(t_1,\dots,t_n;t_{n+1},\dots,t_{n+m})$ as
follows: If $B$ is not the only through-block of $\sigma$
then we choose the unique cyclic order on $B$ ({\it c.f.}
section 2.2) and, writing it as a cycle
$B=(i(1),\dots,i(k))$, we put
\begin{equation}\label{csi1}
\csi_B(t_1,\dots,t_n;t_{n+1},\dots,t_{n+m})
:=\csi_1\bigl(t_{i(1)}, t_{i(2)}, \dots, t_{i(k)}\bigr).
\end{equation}
If $B$ is the only through-block of $\sigma$, then we write
it as $B=B_1\cup B_2$ with $B_1=(i(1),\dots,i(k))\subset
[1,n]$ and $B_2=(j(1),\dots,j(l))\subset [n+1,n+m]$, where
we induce the cyclic order of $[1,n]$ on $B_1$ and the
cyclic order of $[n+1,n+m]$ on $B_2$. For such a block $B$
we put
\begin{equation}\label{csi2}
\csi_B(t_1,\dots,t_n;t_{n+1},\dots,t_{n+m})
:=\csi_2\bigl(t_{i(1)},\dots, t_{i(k)};
t_{j(1)},\dots,t_{j(l)}\bigr).
\end{equation}
Finally, we define
\begin{equation}
\csi_\sigma(t_1,\dots,t_n;t_{n+1},\dots,t_{n+m})
:=\prod_{B\in \sigma}\csi_B(t_1,\dots,t_n;t_{n+1},\dots,t_{n+m}).
\end{equation} 
\end{notation}

\noindent
Here are some examples of this notation: Consider $n=m=2$ and 
$$\sigma_1=\{(1,3),(2,4)\}\qquad 
\sigma_2=\{(1,2,3),(4)\}.$$
Then $\sigma_1$ has two through-blocks so
\begin{equation}
\csi_{\sigma_1}(t_1,t_2;t_3,t_4) = \csi_1(t_1,t_3)\csi_1(t_2,t_4)
\end{equation}
whereas  $\sigma_2$ has one through-block so
\begin{equation}
\csi_{\sigma_2}(t_1,t_2; t_3,t_4) = \csi_1(t_4)
\csi_2(t_1,t_2;t_3).
\end{equation}

\begin{theorem}\label{main}

Let $(\cA_N,\ff_N,\nappa_N)$ ($N\in\NN$) be second order
non-com\-mut\-at\-ive probability spaces and let, for each
$N\in\NN$, unital subalgebras
$\cA_N^1,\dots,\cA_N^N\subset\cA_N$ be given which are free
with respect to $(\ff_N,\nappa_N)$.  Let $\cI$ be an index
set and assume that we have, for each $t\in\cI$ and each
$N\in\NN$, elements $$q_N^i(t)\in\cA_N^i\qquad
(i=1,\dots,N),$$ such that the following properties are
satisfied:

\noindent
(a) The distribution of the $q_N^i(t)$ under 
$(\ff_N,\nappa_N)$ is invariant under permutations of the upper indices, 
i.e., for all $N\in\NN$, and all permutations
$\pi:[1,N]\to[1,N]$ we have for all $n,m\in\NN$, $t_1,\dots,t_{n+m}\in
\cI$ and all $i(1),\dots,i(n+m)\in [1,N]$ that
$$
\ff_N\bigl(q_N^{i(1)}(t_1)\cdots q_N^{i(n)}(t_n)\bigr)=
\ff_N\bigl(q_N^{\pi\circ i(1)}(t_1)\cdots q_N^{\pi\circ
i(n)}(t_n)\bigr)
$$
and
\begin{multline*}
\nappa_N
\bigl(q_N^{i(1)}(t_1)\cdots q_N^{i(n)}(t_n),q_N^{i(n+1)}(t_{n+1})\cdots
q_N^{i(n+m)}(t_{n+m})\bigr)\\=
\nappa_N\bigl(q_N^{\pi\circ i(1)}(t_1)\cdots q_N^{\pi\circ
i(n)}(t_n),q_N^{\pi\circ i(n+1)}(t_{n+1})\cdots
q_N^{\pi\circ i(n+m)}(t_{n+m})\bigr)
\end{multline*}
(b) For all $n,m\in\NN$ and 
all $t_1,\dots,t_n,t_{n+1},\dots,t_{n+m}\in\cI$ 
there exist constants $\csi_1(t_1,\dots,t_n)$
and $\csi_2(t_1,\dots,t_n;t_{n+1},\dots,t_{n+m})$ such that
\begin{equation}\label{assumpt1}
\lim_{N\to\infty}N\cdot \ff_N\bigl(q_N^i(t_1)\cdots q_N^i(t_n)\bigr)= 
\csi_1(t_1,\dots,t_n)
\end{equation}
and
\begin{multline}\label{assumpt2}
\lim_{N\to\infty} N\cdot
\nappa_N\bigl
(q_N^i(t_1)\cdots q_N^i(t_n),q_N^i(t_{n+1})\cdots q_N^i(t_{n+m})\bigr)\\=
\csi_2(t_1,\dots,t_n;t_{n+1},\dots,t_{n+m}).
\end{multline}
For $t\in\cI$ and $N\in\NN$ let
$$S_N(t):=q_N^1(t)+\cdots+ q_N^N(t)\in\cA_N.$$
Then we have
\begin{multline}\label{limit}
\lim_{N\to\infty}\nappa_N\bigl(S_N(t_1)\cdots S_N(t_n),S_N(t_{n+1})\cdots
S_N(t_{n+m})\bigr)\\
=\sum_{\sigma\in NC(n,m)}\csi_\sigma(t_1,\dots,t_n;t_{n+1},
\dots,t_{n+m})
\end{multline}
\end{theorem}

Note that the left-hand side of the expressions
(\ref{assumpt1}) and (\ref{assumpt2}) are independent of the
value of the index $i$, and that the functions $\csi_1$ and
$\csi_2$ defined there have the traciality properties which
are required in Notation \ref{not:oje}.

\begin{proof}
For better legibility, we will suppress in the following the
index $N$ at $\ff_N$ and $\nappa_N$ and just write $\ff$ and
$\nappa$, respectively.

We have
\begin{align*}
&\nappa\bigl(S_N(t_1)\cdots S_N(t_n),S_N(t_{n+1})\cdots
S_N(t_{n+m})\bigr)\\
&=\sum_{i:[1,n+m]\to[1,N]}
\nappa\bigl(q_N^{i(1)}(t_1)\cdots q_N^{i(n)}(t_n),
q_N^{i(n+1)}(t_{n+1})\cdots q_N^{i(n+m)}(t_{n+m})\bigr)
\end{align*}
Because of our invariance assumption (a), the value of the
term
\begin{equation}\label{kappa-value}
\nappa\bigl(q_N^{i(1)}(t_1)\cdots q_N^{i(n)}(t_n),
q_N^{i(n+1)}(t_{n+1})\cdots q_N^{i(n+m)}(t_{n+m})\bigr)
\end{equation}
depends on $i$ only through the information where these
indices are the same and where they are different. As usual, this
information is encoded in a partition $\sigma$ of the set $[1,n+m]$, and 
we denote the common value of (\ref{kappa-value}) for all $i$ with
$\ker(i)=\sigma$ by
\begin{equation}\label{kappa-part}
\nappa_\sigma\bigl(q_N(t_1)\cdots q_N(t_n),
q_N(t_{n+1})\cdots q_N(t_{n+m})\bigr).
\end{equation}
Then we can continue our calculation as follows:
\begin{align*}
&\nappa\bigl(S_N(t_1)\cdots S_N(t_n),S_N(t_{n+1})\cdots
S_N(t_{n+m})\bigr)\\
&=\sum_{\sigma\in\cP(n+m)}\mathop{\sum_{i:[1,n+m]\to[1,N]}}_{\ker(i)=\sigma}
\nappa\bigl(q_N^{i(1)}(t_1)\cdots q_N^{i(n)}(t_n), \\[-10pt]
&\kern140pt q_N^{i(n+1)}(t_{n+1})\cdots q_N^{i(n+m)}(t_{n+m})\bigr)\\
&=\sum_{\sigma\in\cP(n+m)}
\nappa_\sigma\bigl(q_N(t_1)\cdots q_N(t_n),
q_N(t_{n+1})\cdots q_N(t_{n+m})\bigr)\cdot (N)_{\vert\sigma\vert},
\end{align*}
because the number of $i:[1,n+m]\to[1,N]$ with the property
$\ker(i)=\sigma$ is given by
$$N(N-1)\cdots (N-\vert\sigma\vert+1)=: (N)_{\vert\sigma\vert}.$$

We have now to examine the contributions for different
$\sigma$.  Let us first assume that $\sigma$ contains a
block $B$ which is either contained in $[1,n]$ or contained
in $[n+1,\dots,n+m]$ and all of whose elements are
consecutive in the induced cyclic order.  Because of
traciality of $\nappa$ it suffices to consider the case
$B=[1,s]$ for some $s$ with $1\leq s\leq n$.  By Lemma
\ref{properties}, this implies
$$\displaylines{
\nappa_\sigma\bigl(q_N(t_1)\cdots q_N(t_s)\cdots q_N(t_n),
q_N(t_{n+1})\cdots q_N(t_{n+m})\bigr) \cr
=\ff\bigl(q_N(t_1)\cdots q_N(t_s)\bigr)\cdot
\nappa_{\sigma'}\bigl(q_N(t_{s+1})\cdots q_N(t_n),
q_N(t_{n+1})\cdots q_N(t_{n+m})\bigr),
\cr}$$
where $\sigma'$ is that partition which results from
$\sigma$ by removing the block $B=[1,s]$ and relabelling
elements.  Since
$$
\lim_{N\to\infty}N\cdot \ff(q_N(t_1)\cdots q_N(t_s))= 
\csi_1(t_1,\dots,t_s),
$$
the block $B$ makes exactly the contribution to the final
result as claimed in Eq. (\ref{limit}). Thus, by
successively removing such blocks, it suffices to consider
$\sigma$'s which have no
blocks which are contained in either $[1,n]$ or
$[n+1,n+m]$ and which consist of cyclically consecutive
elements. 

So let us now assume that $\sigma$ contains 
no blocks which are contained in either $[1,n]$ or
$[n+1,n+m]$ and which consist of cyclically consecutive
elements, and consider (\ref{kappa-part}).  By
multiplying together neighbouring elements corresponding to
the same block of $\sigma$ we can rewrite the two arguments
of $\nappa$ in a cyclically alternating form. The fact that
$\sigma$ contains no blocks of the form treated above
implies that after this rewriting of arguments each involved
subalgebra occurs at least twice. But then Lemma
\ref{fact:bad} implies that, unless condition (*) is
satisfied, we can write all these terms as sums of products
of at least $\vert\sigma\vert +1$ factors. By our
assumption, each of these factors multiplied by $N$
converges to a finite number; however, since we have more
than $\vert\sigma\vert$ factors, this product multiplied by
$N^{\vert\sigma\vert}$ will vanish in the limit
$N\to\infty$. This means that we can only get a
non-vanishing limit for a $\sigma$ which satisfies condition
(*) of Lemma \ref{fact:bad}. However, these are exactly the
cases where each block $B$ of $\sigma$ is of the form
$B=B_1\cup B_2$, where $B_1\subset [1,n]$ and
$B_2\subset[n+1,m+1]$, are non-empty, and each consists of
consecutive numbers with respect to the inherited order.
Furthermore, the cyclic order of the restrictions of all
blocks to the interval $[1,n]$ must be the inversion of the
cyclic order of the restrictions of all blocks to the
interval $[n+1,n+m]$.  In this case (\ref{kappa-part})
calculates as follows. If we have only one block in
$\sigma$, then our assumption, Equation (\ref{assumpt2}),
gives, in the limit, for such a $\sigma$ the contribution
$$\csi_2(t_1,\dots,t_n;t_{n+1},\dots,t_{n+m}).$$
If, on the other side, $\sigma$ has more than one block,
then we get, according to the description of annular non-crossing
 partitions in section 2.2 and our assumption
(\ref{assumpt2}), the product of $\csi_B$ over all blocks
$B$ of $\sigma$.

Note that the reduction above leads to non-vanishing
contributions exactly for non-crossing partitions $\sigma$
from $NC(n,m)$ and each such partition $\sigma$ contributes
a term $\csi_\sigma(t_1,\dots,t_n;t_{n+1},\dots,t_{n+m})$.
\end{proof}

\section{Proofs of Theorems \ref{main1} and \ref{main2}}

Now we can prove our main theorems by reducing them to the
situation covered in our limit theorem.

\subsection{Proof of Theorem \ref{main1}}

We have to show that for all $n,m\in\NN$ and
$f_1,\dots,f_{n+m}\in\HH_\RR$
$$\displaylines{
\la \cyc \omega(f_1) \cdots \omega(f_n)\Omega,
\cyc \omega(f_{n+m})\cdots \omega(f_{n+1})\Omega\ra_{\cyclic} \hfil\cr
\mbox{} =
\sum_{\pi\in NC_2(n,m)}\prod_{(i,j)\in\pi}\la f_i,f_j\ra.
\cr}$$

Note that we can, for any $N\in\NN$, replace $\HH$ by
$$\bigoplus_{i=1}^N \HH=
\HH\oplus\dots\oplus \HH \qquad \text{($N$ summands)}$$ 
and $\omega(f)$ by
$$\frac 1{\sqrt{N}} \omega(f\oplus\dots\oplus f).$$

We can then put this into the framework of our general limit
theorem by letting $\cI=\HH_\RR$,
$$\cA_N=\cA(\bigoplus_{i=1}^N \HH),\qquad
\cA_N^i=\cA(0\oplus\cdots\oplus
\underset{\text{$i$-th}}{\HH}\oplus\cdots\oplus 0)
$$
$$\ff_N(a)=\la a\Omega,\Omega\ra\qquad
(a\in \cA_N)$$
and
$$\nappa_N(a,b)=\la \cyc a\Omega,\cyc b^*\Omega\ra_{\cyclic}
\qquad (a,b\in\cA_N)$$
and finally, for $f\in\HH_\RR$,
$$q_N^i(f)=\frac1{\sqrt N} \omega(0\oplus\cdots\oplus
\underset{\text{$i$-th}}{f}\oplus\cdots\oplus 0)\in\cA_N^i.$$
Let us check that $\cA_N^1,\dots,
\cA_N^N\subset\cA_N$ are free with respect to $(\ff_N,\nappa_N)$:
Freeness with respect to $\ff_N$ is well-known, so we only have
to consider $\nappa_N$.
Take centered and cyclically alternating tuples $(a_1,\dots,a_n)$ and
$(b_1,\dots,b_m)$ from $\cA_N$. 
Let us only consider the case $n,m\geq 2$, the cases were at least
one of them is 1 are similar.
Note that the centeredness of the
$a_i$ implies that each $a_i\Omega$ has no component in the direction
$\Omega$ and thus, by the fact that neighbours are from algebras
with orthogonal Hilbert spaces, we have 
$$a_1a_2\cdots a_n\Omega=(a_1\Omega)\otimes (a_2\Omega)\odo
(a_n\Omega).$$ Since also the first and the last element are
orthogonal, the action of $\cyc$ becomes in this case just
$$\cyc a_1 a_2\cdots a_n\Omega=\bigl[a_1\Omega\otimes a_2\Omega
\odo a_n\Omega\bigr].$$
In the same way we have
$$\cyc b_m^\ast b_{m-1}^\ast \cdots b_1^\ast
\Omega=\bigl[b_m^\ast \Omega\otimes b_{m-1}^\ast \Omega\odo
b_1^\ast\Omega\bigr].$$
If we take now the inner product in the cyclic Fock space between
these two vectors, then we get
\begin{align*}
\nappa_N(a_1\cdots a_n,b_1\cdots b_m)&=
\delta_{nm}\sum_{k=0}^{n-1}\la a_1\Omega,b_{n+k}^\ast
\Omega\ra\cdots \la a_n\Omega,b_{1+k}^\ast\Omega\ra\\
&=\delta_{nm}\sum_{k=0}^{n-1} \ff(a_1b_{n+k})\cdots
\ff(a_nb_{1+k}),
\end{align*}
as required by our Definition \ref{free}. Thus the
subalgebras $\cA_N^1,\dots,\cA_N^N$ are free with respect to
$(\ff_N,\nappa_N)$.  The invariance assumption on the
distribution with respect to $(\ff_N,\nappa_N)$ is also
easily verified and so we can apply our limit theorem.

Let $S_N(f)=q_N^1(f)+\cdots+ q_N^N(f)$. Since
$$\displaylines{
\langle \cyc \omega(f_1) \cdots \omega(f_n) \Omega, 
        \cyc\omega(f_{n+m}) \cdots \omega(f_{n+1}) \rangle_{\cyclic} \hfill\cr
\mbox{} = 
\nappa_N(S_N(f_1) \cdots S_N(f_n), S_N(f_{n+1}) \cdots S_N(f_{n+m}))
\cr}$$
we can take the limit as $N \rightarrow \infty$ and apply
Theorem \ref{main}. So it remains to identify the limits
$\csi_1$ and $\csi_2$ in the assumption of that theorem.
One sees easily that
\begin{eqnarray*}
\csi_1(f_1,\dots,f_n) &=&
\lim_{N\to\infty} N\cdot \ff_N\bigl(q_N^i(f_1)\cdots q_N^i(f_n)\bigr) \\
\mbox{} &=&
\begin{cases}
\la f_1,f_2\ra & \text{if $n=2$}\\
0 & \text{otherwise}
\end{cases}
\end{eqnarray*}
and
\begin{align*}
\csi_2(f_1,\dots, & f_n;g_1,\dots,g_m) \\
&=
\lim_{N\to\infty} N\cdot
\nappa_N\bigl(q_N^i(f_1)\cdots q_N^i(f_n),q_N^i(g_1)\cdots q_N^i(g_m)
\bigr)\\&=
\begin{cases}
\la f_1,g_1\ra & \text{if $n=1=m$}\\
0 & \text{otherwise}
\end{cases}.
\end{align*}
This gives exactly our claim. \qed

\subsection{Proof of Theorem \ref{main2}}

From equation (\ref{wishartcumulant}) we only have to prove that
\begin{multline}
\la\cyc p(d_1)\cdots p(d_n)\Omega,\cyc p(d_{n+m}^\ast) \cdots
p(d_{n+1}^\ast)\Omega\ra_{\cyclic}\\
=\sum_{\pi\in \SNC(n,m)}\psi_\pi(d_1,\dots,d_n, d_{n+1},\dots,d_{n+m}).
\end{multline}
Note that we can replace $\cC$ by $\cC\otimes L^\infty[0,1]$,
$\psi$ by $\psi \otimes \tau$, where $\tau$ is integration with respect to
Lebesgue measure on $[0, 1]$, and
for each $N\in\NN$, $p(d)$ by
$$
p(d\otimes\chi(0,1))=p_N^1(d)+p_N^2(d)+\dots+p_N^N(d),
$$
where we have put
$$p_N^i(d):=p(d\otimes \chi(I_N^i))$$
with $\chi(I_N^i)$ denoting the characteristic function of the interval
$$I_N^i= \Big(\frac{i-1}N,\frac iN \Big).$$
This fits into the framework of our general limit theorem by
putting $\cI=\cC$,
$$\cA_N=\cA(\cC\otimes L^\infty(0,1)),\qquad
\cA_N^i=\cA(\cC\otimes L^\infty(I_N^i)),$$
$$\ff(a)=\la a\Omega,\Omega\ra,
\qquad
\nappa(a,b)=\la \cyc a\Omega, \cyc b^*\Omega\ra_{\cyclic}
\qquad(a,b\in\cA_N),$$
and finally
$$q_N^i(d)=p_N^i(d)\in\cA_N^i.$$
One can check again by the same arguments as for the
semi-circular case that $\cA_N^1,\dots,\cA_N^N\subset\cA_N$
are free with respect to $(\ff_N,\nappa_N)$. Also the
invariance assumption on the distribution with respect to
$(\ff_N,\nappa_N)$ is easily verified.

Since $p(d)$ has, for each $N$, the same moments with
respect to $\ff_N$ and $\nappa_N$ as
$S_N(d)=q_N^1(d)+\cdots+q_N^N(d)$, we can calculate the
moments of $p(d)$ via $S_N(d)$ by sending $N\to\infty$ and
invoking our limit theorem, Theorem \ref{main}. It only
remains to identify the limits $\csi_1$ and $\csi_2$ from
the hypothesis of the theorem, and show that
\begin{eqnarray}\label{e41}\nonumber\lefteqn{
\sum_{\sigma \in NC(n,m)} \csi_\sigma(d_1,\dots,d_n; d_{n+1},\dots,
d_{n+m})  }\\ 
&=& \sum_{\pi\in \SNC(n,m)}\psi_\pi(d_1,\dots,d_n,
d_{n+1},\dots,d_{n+m}) 
\end{eqnarray}

Note that each inner product appearing in the calculation of
$$\ff_N(p_N^i(d_1) \cdots p_N^i(d_n))$$
gives a factor $1/N$; one inner product must be involved in
any case to get a non-vanishing result, thus the sought
limits single out exactly the contributions with one inner
product.  In the case of $\ff_N(p_N^i(d_1)\cdots \ab
p_N^i(d_n))$ this means that $p_N^i(d_n)$ must act as a
creation operator, $p_N^i(d_1)$ as an annihilation operator
and all the other $p$'s as preservation operators, thus
\begin{align*}
\csi_1(d_1, \dots , d_n) &=
\lim_{N\to\infty} N\cdot\ff_N(p_N^i(d_1)\cdots p_N^i(d_n)) \\
&= \la d_2d_3\cdots d_n,d_1^*\ra =\psi(d_1d_2\cdots d_n)
\end{align*} 
In the case of $\nappa_N$ one has to note that the only
relevant contributions to $\cyc p_N^i(d_1)\cdots
p_N^i(d_n)\Omega$ are of the form: $p_N^i(d_n)$ must act as
creation operator; since $\cyc \Omega=0$, no annihilation
operator can appear, but since $\cyc$ can also act by
multiplication of arguments there might be a second action
as creation operator (let's say of $p_N^i(d_k)$), all the
other $p$ have to act as preservation operators. Thus the
relevant contributions of $\cyc p_N^i(d_1)\cdots
p_N^i(d_n)\Omega$ are the terms with $k=1,\dots,n$ of the
form
$$\cyc (d_1\cdots d_k\otimes d_{k+1}\cdots d_n).$$
Since we are looking for terms which give in the end exactly one
inner product, the relevant action of $\cyc$ 
is given by multiplying arguments and yields terms of the form 
$$[d_{k+1}\cdots d_n d_1\cdots d_k]\qquad\text{for
some $k=1,\dots,n$.}$$
In the same way the relevant contributions of 
$\cyc p_N^i(d_{n+m}^\ast)\cdots p_N^i(d_{n+1}^\ast)\Omega$ are of the form
$$[d_{n+l-1}^\ast \cdots d_{n+1}^\ast d_{n+m}^\ast d_{n+m-1}^\ast \cdots
d_{n+l}^\ast]\qquad \text{for some $l=1,\dots,m$.}$$ 

\noindent
Thus we have 
\begin{align}\label{csi3}
\csi_2&(d_1, \dots , d_n; d_{n+1}, \dots , d_{n +
m}) \\ \nonumber
= &
\lim_{N\to\infty} N\cdot\nappa_N(p_N^i(d_1)\cdots p_N^i(d_n),
p_N^i(d_{n+1})\cdots p_N^i(d_{n+m}))  \\ \nonumber
=&
\lim_{N\to\infty} N\cdot \la 
\cyc p_N^i(d_1)\cdots p_N^i(d_n)\Omega,
\cyc p_N^i(d_{n+m}^\ast)\cdots
p_N^i(d_{n+1}^\ast)\Omega\ra_{\cyclic}\\ \nonumber
=&\sum_{k=1}^n\sum_{l=1}^m \la  [d_{k+1}\cdots d_n d_1\cdots d_k],
[d_{n+l-1}^\ast \cdots d_{n+1}^\ast d_{n+m}^\ast \cdots
d_{n+l}^\ast]\ra_\cyclic\\ \nonumber
=&\sum_{k=1}^n\sum_{l=1}^m
\psi(d_{n+l}\cdots d_{n+m} d_{n+1}\cdots  d_{n+l-1}d_{k+1}\cdots
d_nd_1\cdots d_k)  \end{align}

Suppose $\sigma \in NC(n, m)$ has more than one through-block. Then 
for each block $B$, $\csi_B (d_1, \dots , d_n; d_{n+1}, \dots ,
d_{n+m}) = \psi_B( d_1, \dots , d_{n+m})$ by equation
(\ref{csi1}). Thus
$$
\csi_\sigma(d_1,\ab \dots d_n; d_{n+1}, \dots , d_{n+m}) = 
\psi_\pi(d_1, \dots d_n, d_{n+1}, \dots , d_{n+m})
$$
where $\pi \in S_{NC}(n, m)$ is the unique permutation whose cycle
decomposition is the partition $\pi$.

\noindent
\BoxedEPSF{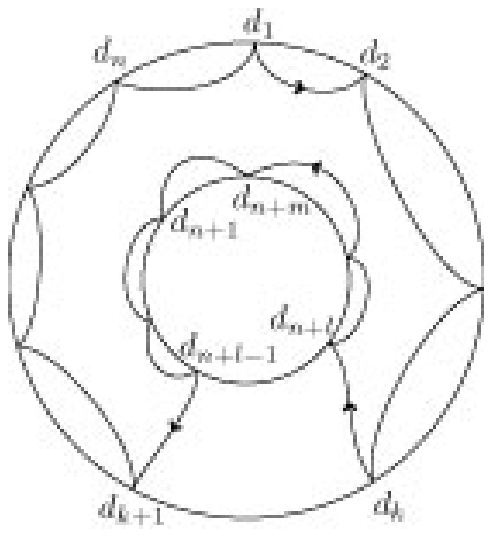}\hfil $\vcenter to
150pt{\hsize150pt\raggedright\noindent\vfill\small {\bf
Figure \figno.} One of the possible cycles which contribute
to the sum in $\csi_B$ in the case $B$ is a unique through
block.\vfill}$

\medskip

Now suppose that $\sigma$ has only one through-block, let
$[\sigma]$ be the set of all $\pi \in S_{NC}(n, m)$ whose
cycle decomposition gives the partition $\sigma$. If $B$ is
a block of $\sigma$ which is not a through-block then again
by equation (\ref{csi1}) $\csi_B$ and $\psi_B$ are equal. If
$B$ is the unique through-block then as in equation
(\ref{csi2}) write $B = \{j_1, \dots , j_r\} \cup \{j_{r+1},
\dots , j_{r+s} \}$. Then by (\ref{csi3}) $$\displaylines{
\csi_B(d_{j_1}, \dots , d_{j_r}; d_{j_{r+1}}, \dots ,
d_{j_{r+s}}) \mbox{} = \csi_2( d_{j_1}, \dots , d_{j_r};
d_{j_{r+1}}, \dots , d_{j_{r+s}}) \cr \mbox{} = \sum_c
\psi_c( d_{j_1}, \dots, d_{j_r}, d_{j_{r+1}}, \dots,
d_{j_{r+s}} ) \cr}$$ where $c$ runs over the cycles in $\pi
\in [\sigma]$ which give the block $B$. Hence $\csi_\sigma =
\sum_{\pi \in [\sigma]} \psi_\pi$ and thus equation
(\ref{e41}) is proved. \qed

\section{Diagonalization of fluctuations}

Let us now use our description of fluctuations of random
matrices in terms of operators to diagonalize these
fluctuations.  The one-dimensional Gaussian case is well
established in the physical and mathematical literature
(see, e.g., \cite{Pol,AJM,Joh}), whereas looking on the
one-dimensional Wishart case and, in particular, on the
multi-dimensional Gaussian case was initiated by
Cabanal-Duvillard \cite{C-D}. Indeed, trying to understand
and reproduce the results of Cabanal-Duvillard was the
original motivation for our investigations.
 
Since the fluctuations are given by taking inner products in
cyclic Fock space, we can achieve such a diagonalization by
taking functions of our operators which yield elementary
tensors in cyclic Fock space. This means we are looking for
a kind of cyclic Wick products.

\subsection{Semi-circular case}

We should look for cyclic analogues of the Wick products
$W(f_1\odo f_n)$.  Let us denote them by $C(f_1 \odo
f_n)$. They should be determined by the property that
$$\cyc C(f_1\odo f_n)\Omega=[f_1 \odo f_n].$$

Notice that we have
\begin{align*}
\cyc W&(f_1\odo f_n)\Omega = \cyc (f_1\otimes\cdots\otimes f_n) \\
&=[f_1\odo f_n]+\la f_1,\bar f_n\ra\cdot\cyc f_2\otimes\cdots\otimes
f_{n-1}\\
&=\cyc C(f_1 \odo f_n)\Omega
+\la f_1,\bar f_n\ra\cdot\cyc W(f_2 \odo f_{n-1})\Omega ,
\end{align*}
thus we could define these cyclic Wick products by the following
recursion:
$$C(f_1 \odo f_n)=W(f_1 \odo f_n)-\la f_1,\bar f_n\ra\cdot
W(f_2,\dots,f_{n-1}).$$
For $n=1$, this means, of course,
$$C(f)=W(f)=\omega(f).$$
If we put
$$f:=f_1=f_2=\dots =f_n\qquad (\text{with $\Vert f\Vert =1$),}$$
then we know that 
$$W(f^{\otimes n})=U_n(\omega(f)/2),$$
where the $\{U_n\}_n$ are the Chebyshev polynomials of the second
kind. Let $V_n(x) = U_n(x/2)$; then $V_n(\omega(f)) =
W(f^{\otimes^n})$. Now, if we write our cyclic Wick polynomials in
this one-dimensional case as
$$C(f^{\otimes n})= 2\, T_n(\omega(f)/2),$$
then these $T_n$ must satisfy
$$ 2 T_n=U_n-U_{n-2}\qquad (n\geq 2)$$
and 
$$T_1(x) = U_1(x)/2 = x.$$
This shows that the $\{T_n \}$ are Chebyshev polynomials of the first
kind.

Let us now consider the multi-dimensional case. It is easy
to see that if $f_i$ is orthogonal to $f_{i+1}$ for all
$i=1,\dots,k$, then we have for all $n(1),\dots,n(k)>0$ that
$$
W(f_1^{\otimes n(1)}\otimes f_2^{\otimes 
n(2)}\otimes\cdots\otimes f_k^{\otimes n(k)})=
W(f_1^{\otimes n(1)})\cdot W(f_2^{\otimes n(2)})\cdots W(f_k^{\otimes
n(k)}).
$$
If we assume in addition that also $f_1$ and $f_k$ are
orthogonal then we get for the corresponding $C$:
\begin{align*}
C(f_1^{\otimes n(1)}\otimes f_2^{\otimes
n(2)} \odo & f_k^{\otimes n(k)})=
W(f_1^{\otimes n(1)}\otimes f_2^{\otimes 
n(2)}\otimes\cdots\otimes f_k^{\otimes n(k)})\\
&=W(f_1^{\otimes n(1)})\cdot W(f_2^{\otimes n(2)})\cdots W(f_k^{\otimes
n(k)}).
\end{align*}
The covariance between such functions in our random matrices
is given by the inner product in the cyclic Fock space.  If
we have $k,l\geq 2$ and $f_1,\dots,f_k\in\HH_\RR$ and
$g_1,\dots,g_l\in\HH_\RR$ such that $f_i\perp f_{i+1}$ for
$i=1,\dots,k$ and $g_i\perp g_{i+1}$ for $i=1,\dots,l$ then
we have
$$\displaylines{
\lim_{N\to\infty}
\kk_2 \Bigl\{ \Tr[V_{n(1)}(X_N(f_1)) \cdots V_{n(k)}(X_N(f_k)],
\hfil\cr
\hfill\Tr[V_{m(1)}(X_N(g_1)) \cdots V_{m(l)}(X_N(g_l))]\Bigr\} \cr
\mbox{} =
\la [f_1^{\otimes n(1)}\odo f_k^{\otimes n(k)}],
[g_l^{\otimes m(l)}\odo g_1^{\otimes m(1)}]\ra_\cyclic.
\cr}$$ 
Thus we recover the results of Cabanal-Duvillard \cite{C-D}
for that case.

\subsection{Compound Poisson case}

Again, we are looking for polynomials $C( d_1\odo d_n)$
which have the property
$$\cyc C(d_1\odo d_n)\Omega=[d_1\odo d_n].$$
We have
\begin{align*}
&\cyc W(d_1\odo d_n)\Omega=\cyc (d_1 \odo d_n) \\ 
&\mbox{} =
[d_1,\dots,d_n]+[d_nd_1,d_2,\dots,d_{n-1}]
+ \psi(d_1d_n) \, \cyc (d_2 \odo d_{n-1}) \\ 
&\mbox{}=
\cyc C(d_1\odo d_n)\Omega + \cyc C(d_nd_1\otimes d_2\odo d_{n-1})\Omega \\
&\qquad+\psi(d_1d_n)\cyc W(d_2 \odo d_{n-1})\Omega
\end{align*}
Thus we define the $C$'s in the following recursive way:
\begin{multline}
W(d_1\odo d_n) = C(d_1 \odo d_n) \\ 
+ C(d_nd_1\otimes d_2 \odo d_{n-1}) + \psi(d_1d_n) W(d_2\odo d_{n-1})
\end{multline}

There does not seem to be a nice closed form for this in the
one-dimensional case.

Let us also look at the multi-dimensional situation.  We
model this by assuming that we have elements
$d_1,\dots,d_r\in\cC$ such that $d_id_j=0$ for $i\not= j$.
Then we have again for $i(j)\not=i(j+1)$ ($j=1,\dots,n$) and
$k(1),\dots,k(n)>0$ that
$$W(d_{i(1)}^{\otimes k(1)}\odo d_{i(n)}^{\otimes k(n)})=
W(d_{i(1)}^{\otimes k(1)})\cdots W(d_{i(n)}^{\otimes k(n)}).$$
If also $i(1)\not= i(n)$, then we have again equality
between $W$ and $C$, i.e.
$$\displaylines{
C(d_{i(1)}^{ k(1)}\odo d_{i(n)}^{ k(n)})=
W(d_{i(1)}^{\otimes k(1)}\odo d_{i(n)}^{\otimes k(n)}) \cr
\hfill =
W(d_{i(1)}^{\otimes k(1)})\cdots W(d_{i(n)}^{\otimes k(n)}).
\cr}$$

\subsection{Poisson case}

Let us specialize the general compound Poisson 
case to the usual Poisson case.

The usual Poisson case is special within the class of
compound ones by a very special state on $\cC$.  Restrict
for the moment to one random matrix, i.e., the algebra $\cC$
is generated by one element $d$.  Then the fact that we have
a free Poisson variable $p(d)$ means that this $d$ is a
projection and thus
$$\psi(d^k)=\psi(d)=:\lambda.$$
So we can identify $p(d)=p(d^2)=\dots$ and everything
reduces again to polynomials in just one variable $p(d)$.
Again one knows (see \cite[Theorem 4.11]{Ans}) that the linear
Wick polynomials
$W_n(d):=W(d^{\otimes n})$ are given by the orthogonal
polynomials with respect to the distribution of $p(d)$
(i.e. with respect to the Marchenko-Pastur = free Poisson
distribution). Let us denote these polynomials by $\Pi_n$,
then we have
$$W_n(d)=\Pi_n(p(d)).$$
If we put $C_n(d):=C(d^{\otimes n})$, then the general
relation between $W$ and $C$ becomes in this case:
$$W_n(d)=C_n(d)+C_{n-1}(d)+\lambda W_{n-2}.$$
If we put $C_n(d)=\Gamma_n(p(d))$ for some polynomials $\Gamma_n$,
then the above tells us that
$$\Pi_n-\lambda \Pi_{n-2}=\Gamma_n+\Gamma_{n-1}.$$
This gives us exactly the polynomials $\{ \Gamma_n \}$ which
appear in Cabanal-Duvillard's results \cite{C-D}.

As an extension of this, we also get the multi-dimensional
Poisson case: There the ``diagonalizing polynomials" in more
than one variable are given by alternating products in the
one-dimensional linear polynomials $\{ \Pi_n \}$.

A more detailed investigation of this diagonalization of
fluctuations will be presented in \cite{KMS}.

\section{Asymptotic freeness of Gaussian and constant matrices}

Our results about compound Wishart matrices can be
considered as describing the limiting relation between
Gaussian random matrices and constant matrices for special
moments -- namely those with patterns of the form
$X^\ast D_1XX^\ast D_2X\cdots X^\ast D_nX$. This raises, of
course, the question whether we can say something substantial
about the general relation between Gaussian and constant
matrices. In view of the basic theorem of Voiculescu that
Gaussian random matrices and constant matrices are
asymptotically free, we would expect that we should have the
same kind of statement also on the level of fluctuations.  We
want to indicate here that this is indeed the case, thus
providing strong evidence that our notion of ``second order
freeness" is indeed the correct concept.  Note that in the
following definition we make a quite strong requirement on the
vanishing of the higher order cumulants. This is however in
accordance with the observation that in many cases the
unnormalized traces converge to Gaussian random variables.  Of
course, if we have a non-probabilistic ensemble of constant
matrices, then the only requirement is the convergence of
$\kk_1$; all other cumulants are automatically zero.

\begin{definition} \label{def:asymptoticfreeness}
1) Let $\{A_1,\dots,A_s\}_N$ be a sequence of $N\times N$-random
matrices. We say that they have a \emph{second order limit
distribution} if there exists a second order non-commutative
probability space $(\cA,\ff,\ab \nappa)$ and $a_1,\dots,a_s\in\cA$
such that for all polynomials $p_1,p_2,\dots$ in $s$ non-commuting
indeterminates we have
\begin{equation}\label{asymp1}
\lim_{N\to\infty}\kk_1\bigl\{\tr[p_1(A_1,\dots, A_s)]\bigr\}=
\ff\bigl(p_1(a_1,\dots, a_s)\bigr),
\end{equation}

\begin{multline}\label{asymp2}
\lim_{N\to\infty} \kk_2\bigl\{\Tr[p_1(A_1,\dots,A_s)],
\Tr[p_2(A_{1},\dots, A_{s})]\bigr\}=\\
\nappa\bigl(p_1(a_1,\dots,a_s);p_2 (a_{1},\dots, a_{s})\bigr),
\end{multline}
and, for $r\geq 3$,
\begin{equation}\label{asymp3}
\lim_{N\to\infty} \kk_r\bigl\{\Tr[p_1(A_1,\dots, A_s)],\dots,
\Tr[p_r(A_1,\dots, A_s)]\bigr\}=0.
\end{equation}
2) We say that two sequences of $N\times N$-random matrices,
$\{A_1,\dots,A_s\}_N$ and $\{B_1,\dots,B_t\}_N$, are
\emph{asymptotically free of second order} if the sequence
$\{A_1,\dots,A_s,B_1,\dots,B_t\}_N$ has a second order limit
distribution, given by $(\cA,\ff,\nappa)$ and $a_1,\dots,a_s,
b_1,\dots,b_t\in\cA$, and if
the unital algebras
$$\cA_1:=\text{alg}(1,a_1,\dots,a_s)\qquad\text{and}\qquad
\cA_2:=\text{alg}(1,b_1,\dots,b_t)$$
are free with respect to $(\ff,\nappa)$. 
\end{definition}

\begin{remark}
Corollary \ref{gaussiansecondorder} shows that a family $\{ X_N(f)
\}_{f \in \HH_{\RR}}$ of Hermitian Gaussian random matrices has a
second order limit distribution. Theorem \ref{main1} identifies the
limiting distribution in terms of cyclic Fock space, and in the
proof of Theorem \ref{main1} we have in addition shown that the
limiting distribution is free of second order in that if $\KK_1,
\dots , \KK_n \subset \HH$ are orthogonal subspaces and $\cA_i$ is
the algebra generated by $\{ \omega(f) \mid f \in \KK_i \}$ then
$\cA_1, \dots , \cA_n$ are free with respect to $(\ff, \nappa)$
where $\ff(a) = \la a \Omega, \Omega \ra$ and $\nappa(a_1, a_2) =
\la \cyc a_1 \Omega, \cyc a_2^\ast \Omega \ra_{\cyclic}$. Thus we
have shown orthogonal families of Gaussian random matrices are
asymptotically free of second order. 
\end{remark}

\begin{remark}
Corollary \ref{wishartsecondorder} showed that if $\{ X_N \}_N$ is
a sequence of complex Gaussian random matrices and 
$P_N(D_i)=X_N^\ast D_i^{(N)}X_N$ where $\{ D_1^{(N)}, D_2^{(N)},
D_3^{(N)}, \dots D_p^{(N)}\}_N$ is a sequence of $N \times N$ complex
matrices which converges in distribution to $(d_1, d_2, \dots d_p)$
in $( \cC, \psi)$ then the family $\{ P_N(D_i) \}_i$ has a limiting
distribution. Theorem \ref{main2} calculates the limiting
distribution in terms of cyclic Fock space. In the proof of Theorem
\ref{main2} we have shown that the limiting distribution is free of
second order in that if $d_i d_j = 0$ for $i \not = j$ and $\cA_i$ is
the algebra generated by $p(d_i)$ then $\cA_1, \dots , \cA_p$ are
free with respect to $(\ff, \nappa)$ where $\ff(a) = \la a \Omega,
\Omega \ra$ and $\nappa(a_1, a_2) = \la \cyc a_1 \Omega, \cyc
a_2^\ast \Omega \ra_{\cyclic}$. Thus we
have shown orthogonal families of Wishart random matrices are
asymptotically free of second order. 

\end{remark}

Now we can address the question of the relation between
Gaussian random matrices and constant matrices. We can even
be more general for the latter and consider random matrices
which are independent from the Gaussian ones.

Let, as usual, $X_N(f)$ ($f\in\HH_\RR$) be a family of
Hermitian Gaussian random matrices
$$X_N(f)=\bigl(x_{ij}(f)\bigr)_{i,j=1}^N,$$ 
as in section \ref{gaussianrm}

\begin{theorem}\label{main:free}

Let $\{X_N(f)\mid f\in\HH_\RR\}_N$ be a sequence of
Hermitian Gaussian $N\times N$-random matrices and
$\{A_1,\dots,A_s\}_N$ a sequence of $N\times N$-random
matrices which has a second order limit distribution. If
$\{X_N(f)\mid f\in\HH_\RR\}_N$ and $\{A_1,\dots,A_s\}_N$ are
independent, then they are asymptotically free of second
order.
\end{theorem}

The proof of this theorem relies on the same kind of
calculations as, for example, in \cite{MN}. Since we do not
want to go into random matrix calculations here, we defer
more details about this to \cite{KMS}.

If the random matrices $\{A_1,\dots,A_s\}$ are non-random
constant matrices with limiting distribution with respect to
the trace, then all $\kk_r$ vanish identically for $r\geq
2$, thus they have a second order limit distribution, and we
get as a corollary of the above that the asymptotic freeness
between Gaussian random matrices and constant matrices
remains also true on the level of fluctuations, i.e., with
respect to our concept of second order freeness.

A more systematic investigation of this concept will be
pursued in forthcoming publications. In particular,
fluctuations of  Haar distributed unitary random  matrices from this
point of view will be treated in \cite{MSS}.

%


\begin{thebibliography}{ABC}

\bibitem[Ans]{Ans} M. Anshelevich: Appell polynomials and their
relatives, Int. Math. Res. Notices, 2004, no.
65, 3469-3531

\bibitem[AJM]{AJM} J. Ambjorn, J. Jurkiewicz and Y.
Makeenko,  Multiloop correlators for two-dimensional
quantum gravity, Phys. Lett. B 251 (1990), 517 -
524.

\bibitem[C-D]{C-D} T. Cabanal-Duvillard: Fluctuations de la loi
empirique de grandes matrices aleatoires. Ann. Inst. H. Poincare
Probab. Statist. 37 (2001), 373-402.

\bibitem[CC]{CC} M. Capitaine and M. Casalis, Asymptotic Freeness by
Generalized Moments for Gaussian and Wishart Matrices. Application to
Beta matrices. Indiana Univ. Math. J., 53 (2004), 397-431.


\bibitem[Eyn]{Eyn} B. Eynard: Random Matrices, Cours de
Physique  Theorique de Saclay, 2000, $\langle${\ttfamily
http://www-spht.cea.fr/articles/t01/014}$\rangle$.

\bibitem[GSS]{GSS} P. Glockner, M. Sch\"urmann, and R. Speicher:
Realization of free white noises. Arch. Math. 58 (1992), 407-416.

\bibitem[GLM]{GLM} P. Graczyk, G. Letac, H. Massam, The Complex Wishart
Distribution and the Symmetric Group, Ann. of Statistics, 31, (2003),
287-309.

\bibitem[HP]{HP} F. Hiai and D. Petz: The semicircle law, free
random variables and entropy. AMS, Providence, RI, 2000.

\bibitem[HT]{HT} U. Haagerup and S. Thorbj{\o}rnsen: A new application
of random matrices: Ext($C^*_{\text{red}}(F_2)$) is not a group, Annals of
Math., 162 (2005), 711-775.

\bibitem[J]{J} S. Janson, Gaussian Hilbert Spaces, Cambridge Tracts
in Mathematics, vol. 129, Cambridge University Press, Cambridge,
1997.

\bibitem[Joh]{Joh} K Johansson: On fluctuations of eigenvalues of
random Hermitian matrices. Duke Math. J. 91 (1998), 151-204.

\bibitem[KMS]{KMS} T. Kusalik, J. Mingo, and R. Speicher,
Orthogonal Polynomials and Fluctuations of Random
Matrices, preprint, math.OA/0503169. 

\bibitem[MN]{MN} J. Mingo and A. Nica: Annular non-crossing
permutations and partitions, and second-order asymptotics for
random matrices, Inter. Math. Res. Notices, 28
(2004), 1413 - 1460.

\bibitem[MSS]{MSS} J. Mingo, P. \'Sniady, and R. Speicher:
Second Order Freeness and Fluctuations of Random Matrices: 
II. Unitary Random Matrices, preprint, math.OA/0405258. 

\bibitem[NSp]{NSp} A. Nica and R. Speicher: Lectures on
the Combinatorics of Free Probability Theory, Paris, 1999

\bibitem[Pol]{Pol} H. D. Politzer: Random matrix description of the
distribution of mesoscopic conductance, Phy. Rev. B, 40 (1989),
11917 - 11919.

\bibitem[Sp1]{Sp1} R. Speicher: A new example of ``Independence"
and ``White Noise". Probab. Th. Rel. Fields 84 (1990), 141-159.

\bibitem[Sp2]{Sp2} R. Speicher: Combinatorial theory of the
free product with amalgamation and operator-valued free
probability theory. Memoirs of the AMS 132 (627), 1998

\bibitem[Voi1]{Voi1} D. Voiculescu: Limit laws for random matrices
and free products, Invent. Math. 104 (1991), 201 - 220.

\bibitem[Voi2]{Voi2} D. Voiculescu: Lectures on free probability
theory. Lecture Notes in Math. 1738 (Springer, 2000), 279-349.

\bibitem[VDN]{VDN} D. Voiculescu, K. Dykema, and A. Nica:
Free Probability Theory, Providence, RI, Amer. Math.
Soc., 1991.

\bibitem[Zvo]{Zvo} A. Zvonkin: Matrix integrals and map enumeration:
an accessible introduction. Math. Comput. Modelling 26 (1997),
281-304.

\end{thebibliography}
\end{document}